\documentclass[10pt]{article}
\usepackage{latexsym,amsmath,amssymb,theorem}
 \usepackage[applemac]{inputenc}  
\usepackage[all]{xy}

\usepackage{color} 
\usepackage{hyperref}
\newenvironment{Proof}{\noindent \bf Proof   \rm}{\hspace*{\fill}
$\square$ \medskip}

{\theorembodyfont{\itshape}
\newtheorem{Lemma}{Lemma}
\newtheorem{Proposition}[Lemma]{Proposition}
\newtheorem{Fact}[Lemma]{Fact}
\newtheorem{Facts}[Lemma]{Facts}

\newtheorem{Theorem}[Lemma]{Theorem}}
{\theorembodyfont{\upshape}

\newtheorem{Definition}[Lemma]{Definition}
\newtheorem{Remark}[Lemma]{Remark}
\newtheorem{Remarks}[Lemma]{Remarks}

\newtheorem{Example}[Lemma]{Example}  
\newtheorem{Examples}[Lemma]{Examples}} 

\newcommand{\tr }{\textcolor{red}}

\newcommand{\ot}{\otimes}

\newcommand{\TV}{\T_\BV}
\newcommand{\TW}{\T_\BW} 
\newcommand{\N}{\mathbb{N}}
\newcommand{\Z}{\mathbb{Z}}
\newcommand{\cn}{\mathcal{N}}
\newcommand{\T}{\mathcal{T}}
\newcommand{\El}{\mathsf{El}}

\newcommand{\BC}{\mathcal{C}}

\newcommand{\BD}{\mathcal{D}}
\newcommand{\B}{\mathcal{B}}

\newcommand{\BB}{\mathcal{B}}
\newcommand{\BA}{\mathcal{A}}

\newcommand{\BR}{\mathcal{R}}
\newcommand{\BS}{\mathcal{S}}
\newcommand{\BT}{\mathcal{T}}
\newcommand{\BU}{\mathcal{U}}
\newcommand{\BE}{\mathfrak{E}}

\newcommand{\Law}{\mathit{Lawv}}

\newcommand{\Adj}{\mathsf{Adj}}
\newcommand{\Ladj}{\mathsf{Ladj}}
\newcommand{\Radj}{\mathsf{Radj}}

\newcommand{\BV}{\mathcal{V}}
\newcommand{\BW}{\mathcal{W}}
\newcommand{\Ring}{{\mathit{Ring}}}

\newcommand{\Mod}{\mathit{Mod}}

\newcommand{\op}{\mathsf{op}}

\newcommand{\ob}{\mathsf{ob}}

\newcommand{\id}{\mathsf{id}}

\newcommand{\ma}{\mathsf{A}}
\newcommand{\mb}{\mathsf{B}}
\newcommand{\mc}{\mathsf{C}}

\newcommand{\Set}{\mathit{Set}}

\newcommand{\Hom}{\mathrm{Hom}}

\newcommand{\End}{\mathrm{End}}

\newcommand{\Alg}{\mathsf{Alg}}
\newcommand{\Ab}{\mathit{Ab}}

\newcommand{\Mon}{\mathit{Mon}}

\newcommand{\Grp}{\mathit{Grp}}
\newcommand{\Lan}{\mathit{Lan}}
\newcommand{\Vect }{\mathit{Vect}_k}
\newcommand{\Bialg}{\mathit{Bialg}}

\newcommand{\Coalg}{\mathsf{Coalg}}
\newcommand{\Hopf}{\mathit{Hopf}}

\newcommand{\lra}{\rightarrow}

\newcommand{\eq }{\mathrm{Eq}\,}
 \newcommand{\Equi }{\mathsf{Equ}}

\newcommand{\xra}{\xrightarrow}

\begin{document}

\title{Internal Coalgebras in Cocomplete Categories:\\
\Large Generalizing the Eilenberg-Watts-Theorem}
  \author{ Laurent Poinsot\\
\small LIPN, UMR CNRS 7030, University   Sorbonne Paris North, \\ \small and CREA, \'Ecole de l'Air, France. \\
{\tt \small laurent.poinsot@lipn.univ-paris13.fr}
  \and  Hans-E. Porst\thanks{Permanent address: Department of Mathematics, University of Bremen, 28359 Bremen, Germany.} \\
\small Department of Mathematical Sciences, University of Stellenbosch, \\ \small Stellenbosch, South Africa. \\
 {\tt\small porst@math.uni-bremen.de} 
 }

\date{}

\maketitle

\begin{abstract} 
The category of internal coalgebras in a cocomplete category $\BC$  with respect to a variety $\BV$ 
is  equivalent to the category of left adjoint functors from $\BV$ into $\BC$. This can be  seen best when considering  such coalgebras as finite coproduct preserving functors  from $\T_\BV^\op$, the dual of  the Lawvere theory of $\BV$, into $\BC$: coalgebras  are  restrictions of  left adjoints and any such  left adjoint is the left Kan extension of a coalgebra along the embedding of $\T_\BV^\op$ into $\Alg\T$.
Since ${_S\Mod}$-coalgebras in the variety ${_R\Mod}$  for  rings $R$ and $S$ are nothing but left $S$-, right $R$-bimodules,  
the  equivalence  above generalizes  the   Eilenberg-Watts Theorem  and all its previous generalizations.
Generalizing and strengthening Bergman's completeness result for categories of internal coalgebras in varieties  we also prove that the category of 
coalgebras in a locally presentable category $\BC$ is locally presentable and comonadic over  $\BC$  and, hence, complete in particular. 
We  show, moreover,  that Freyd's canonical constructions of internal coalgebras in a variety define
left adjoint functors. Special instances of the respective  right adjoints  appear in various algebraic contexts
 and, in the case where $\BV$ is a  commutative variety,  are coreflectors from the category $\Coalg(\BT,\BV)$ into $\BV$.

\end{abstract}

\section*{Introduction}

The  Eilenberg-Watts \cite{Eil,Watts} Theorem states that the only additive and colimit preserving
 functors $L\colon _R\Mod\xra{ }{_S\Mod}$ between the categories of left modules over not necessarily commutative rings $R$ and $S$ are the functors ${M_R}\ot_R -$ for left $S$-, right $R$-bimodules $_SM_R$. By the (dual of) the Special Adjoint Functor Theorem this is equivalent to saying that these tensoring functors are the only additive and left adjoint functors between module categories. Since equivalences between categories  are adjunctions, this result  is of  importance for describing classical Morita theory of rings, one of the early reasons for the interest in it.

Decades later is has been felt necessary to generalize this to additive left adjoint functors $L\colon {_R\Mod} \xra{ }\BC$ for categories $\BC$   more general than module categories. The most recent generalization \cite{N/S} assumes $\BC$   to be a cocomplete additive and $k$-linear category and $R$ the underlying ring of a $k$-algebra (for some commutative ring $k$) in order to prove a non-trivial result about quasi-coherent sheaves. 

The main problem, obviously, then is to replace the concept of bimodule by a more general one. In \cite{BW} and \cite{N/S} this is achieved by choosing left $R$-objects in $\BC$ (see Section \ref{sec:adj} below). The shortcomings of this approach are that this idea neither can be further generalized  to non-additive categories $\BC$ as codomains for $L$ nor to varieties other then module categories as domains of $L$. Note that the  Eilenberg-Watts  Theorem can, alternatively, be read as characterizing the  additive right adjoint functors $K\colon {_S\Mod}\xra{ }{_R\Mod}$ as the functors $\Hom_S(_SM,-)$ for left $S$-, right $R$-bimodules $_SM_R$, and that these the only functors making the following diagram commute:
\begin{equation}\label{Tlift}
\begin{aligned}
\xymatrix@=2.5em{
{ {_S\Mod}}\ar[rr]^{\Hom_S(_SM,-) } \ar[drr]_{_S\Mod(_SM,-)\ \ \ }&& {_R\Mod} \ar@{->}[d]^{|-| }     \\
& &{ \Set}     
}
\end{aligned}
\end{equation}
The functors $\Hom_S(_SM,-)$, thus, are examples  $\BV$-representable functors  in the sense of Freyd \cite{Freyd}, that is,  of functors making the following diagram commute, where $\BC$ is a category with finite coproducts and $|-|$ denotes the forgetful functor of a variety $\BV$.
\begin{equation}\label{T2lift}
\begin{aligned}
\xymatrix@=2.5em{
{ \BC}\ar[r]^{R_\ma } \ar[dr]_{\BC(A,-)}&{\BV} \ar@{->}[d]^{|-| }     \\
 &{ \Set}     
}
\end{aligned}
\end{equation}

Such functors  provide an equivalent way of describing, for any variety $\BV$, internal $\BV$-coalgebras in a category $\BC$ with finite coproducts, that is,  internal $\BV$-algebras in the dual of a category $\BC$ with finite products.
 The latter notion is essentially obtained by translating the concept of equationally defined algebras from the language of universal algebra into that of category theory,  where satisfaction of an equational algebra axiom is expressed by commutativity of a diagram.

Consequently, the $_R\Mod$-coalgebras in $_S\Mod$ are precisely the left $S$-, right $R$-bimodules,  an observation 
Freyd---somewhat surprisingly---did not notice, and this is the appropriate generalization of bimodules needed for generalizing the Eilenberg-Watts theorem. 

As gradually became clear, the concept of internal coalgebras provides a common  perspective on apparently completely unrelated concepts, constructions and results as, for example, the construction of homotopy groups, the change of rings functors in module theory, the Eilenberg-Watts theorem as mentioned above, or the group-like and primitive elements functors of Hopf algebra theory.

Freyd in his paper initiated a systematic study of these structures by what one may call a \em semi-categorical \em approach 
and introduced canonical ways of constructing coalgebras in a variety. Somewhat later, without reference to these papers,  internal coalgebras in the category of commutative rings were introduced and applied in \cite{TW}. More recently gaps in this theory were closed and non trivial results added by Bergmann \cite{Berg1}; most notably, he proved that the category of $\BV$-coalgebras in a variety $\BW$ is complete and its forgetful functor 
into $\BW$ has a right adjoint.  However, his use of algebraic methods only requires a considerable amount of unusual technicalities. So,  from its very beginning, the study of internal coalgebras  appears as an interplay between algebraic and category theoretic concepts and methods.

It is the intention  of this note to show that a purely categorical approach  --- besides providing the appropriate setting for generalizing  the Eilenberg-Watts Theorem ---  simplifies the presentation of the theory of coalgebras and leads to new results, improvements of known ones, and further applications.  For the sake of those readers not too familiar with categorical methods, we include translations of crucial concepts and results into the language of universal algebra. 

We, hence, consider coalgebras in a category $\BC$ with finite coproducts simply as functors $\T^\op\xra{ \ma}\BC$ from the dual $\T^\op$ of a  Lawvere theory into  $\BC$ which preserve finite coproducts. This approach becomes particularly convincing, when $\BC$ not only has finite coproducts but even is cocomplete. Then there results quite obviously an equivalence of coalgebras and left adjoints, which generalizes the Eilenberg-Watts Theorem and can be  expressed  conceptually as follows: 
Up to natural equivalence, $\T$-coalgebras $\T^\op\xra{ \ma}\BC$   are nothing but the restrictions of left adjoints $\Alg\T\xra{L }\BC$ to $\T^\op$, considered as a subcategory of $\Alg\T$, while the left adjoint  $L_\ma$ corresponding to a coalgebra $\ma$ is characterized by any of the following universal properties, where the first one is purely categorical in nature while the second one  is in view of Section \ref{sec:sift} below of an algebraic flavour.   
  \begin{description}
\item[Kan] $L_\ma$ is the left Kan extension of  $\ma$ along the embedding $\T^\op\xra{Y}\Alg\T$.
\item[Sift ]    $L_\ma$ is the  sifted colimit preserving extension of  $\ma$.
\end{description}

Concerning  improvements of the theory of coalgebras we show for example:
\begin{itemize}
\item 
For every Lawvere theory $\T$ the category $\Coalg(\T,\BC)$ of $\T$-coalgebras in a locally presentable category $\BC$ is again a locally presentable category and, hence, complete in particular. This generalizes considerably Bergman's result in \cite{Berg1}. Moreover, its underlying functor into $\Coalg(\T,\BC)\xra{ |\!|-|\!|}\BC$ has a right adjoint (as shown in \cite{Berg1} for $\BC$ being an arbitrary variety, while  the paper \cite{Davis} claims to present a proof   for corings in the category of rings  (see however its review)), and, moreover, 
is comonadic (as shown for $\BC=\Set$ by Isbell \cite{Isbell}). However, not each such comonadic category is of the form $\Coalg(\T,\BC)$.
\item The canonical constructions  of coalgebras  in a variety already mentioned by Freyd are generalized to coalgebras in a locally presentable category. Their  existence in the case of varieties also becomes immediately  evident by our categorical approach, as is the fact that they are  left adjoint functors. The respective right adjoints can be viewed as a generalization of both, the {\em group-like elements} and the {\em primitive elements} functors, known from Hopf algebra theory. 

In the case of a commutative variety $\BV$ with Lawvere theory $\TV$ the canonical  construction provides a full coreflective embedding of $\BV$ into $\Coalg(\TV,\BV)$.  We show by examples that these embeddings may or may not be equivalences. 
\end{itemize}

We also show  that the theory developed can be of use (besides generalizing the Eilenberg-Watts Theorem); we show for example how the following,  familiar algebraic results are obtained easily. 
\begin{itemize}
\item The \em General linear group \em and \em Special linear group \em functors $GL_n, SL_n\colon\Ring\xra{ }\Grp$ have left adjoints.
\item For any  ring homomorphism $\phi\colon R\xra{ } S$ the \em Extension of scalars-\em functor ${_R\Mod}\xra{ }{_S\Mod}$ is left adjoint to the \em Restriction of scalars-\em functor ${_S\Mod}\xra{ }{_R\Mod}$ and this is left adjoint to the \em Coextension of scalars-\em functor ${_R\Mod}\xra{}{_S\Mod}$.
\item Each symmetric  algebra, that is, each free commutative $R$-algebra $\Sigma_\ast(M)$ over an $R$-module $M$, is a Hopf algebra whose set of primitive elements  contains $M$.
\item Moreover we show that our results provide a natural setting for Morita equivalence of arbitrary varieties, generalizing classical Morita theory for module categories.
\end{itemize}

The paper is organized as follows: 

Section \ref{sec:Prel}  is  mainly devoted to  an introduction to the concept of Lawvere theories.  In its part on varieties we improve Freyd's approach to the interpretation of terms in a category with finite products, necessary in this context. 
What seems to be new is the result that the sifted colimit preserving extension of a functor $\T^\op\xra{ }\BC$ is nothing but its left Kan extension along the embedding $\T^\op\hookrightarrow{}\Alg\T$.

Section \ref{sec:IC} deals with the general theory of internal coalgebras. Subsection \ref{sec:def} provides the various possible definitions and proofs of their equivalence, while subsections \ref{sec:prop} and \ref{sec:ex} deal with the properties of the categories of  internal coalgebras and provide a first set of examples.

Section \ref{sec:CV} develops  the theory the categories of coalgebras in varieties. Subsection \ref{sec:co-var} discusses the  canonical construction of coalgebras on free algebras and  its adjunction; it closes with an application of these results  in the context of Kan's construction of cogroups in the category of groups. Subsection \ref{sec:3.2} provides     some non-trivial applications. Subsection \ref{sec:com} deals with coalgebras  in commutative varieties. We show in particular,  that the canonical construction just mentioned, in this case extends to a full coreflective embedding.

\section{Preliminaries}\label{sec:Prel} 

\subsection{Some concepts and notations}\label{not}
We next introduce some concepts and notations not every reader may  be familiar with.
\begin{description}
\item[\sf Semi-additive categories.] A category $\BC$ is called \em semi-additive\em, if it has finite biproducts. Every such category is  enriched over $_c\Mon$. If $\BC$ even is enriched over $\Ab$ it is called \em additive\em.
 \item[\sf Locally presentable categories.] 
Given a regular cardinal $\lambda$, an object $P$ in  a category $\BC$ is called  \em $\lambda$-presentable  \em  iff $\BC(P,-)$ preserves $\lambda$-filtered colimits. A category $\BC$ is called \em locally presentable\em, provided that it is cocomplete and has, for some $\lambda$, a set of $\lambda$-presentable objects of which each $\BC$-objects is a $\lambda$-filtered colimit. Every locally presentable category is complete, well-powered and cowellpowered, and has certain factorization structures of morphisms (see \cite{AR} for details). 

Examples of such categories are all varieties and quasivarieties, the category of partially ordered sets and, more generally, every universal Horn class, and the category of small categories. Also, for  locally presentable categories $\BC$ and $\BD$, the category  $\Ladj(\BD,\BC)$ of left adjoint functors $\BD\xra{ }\BC$  is locally presentable (see \cite{Bird}).

 \item[\sf Functor  coalgebras.]
 The category $\Coalg F$ of \em $F$-algebras in $\BC$ \em for a functor $\BC\xra{ F}\BC$ has as objects all pairs $(C,C\xra{ \gamma}FC)$ and as morphisms $(C,\gamma)\xra{f}(C',\gamma')$ those $\BC$-morphisms $C\xra{ f}C'$ which satisfy $Ff\circ\gamma = \gamma'\circ f$. The obvious forgetful functor $\Coalg F\xra{ }\BC$ reflects isomorphisms and creates colimits and those limits which are preserved by $F$.

\item[\sf Kan extensions.]\label{itemKan} Given functors $\B \xra{ F} \BC$  and $\BB  \xra{ Y} \BD$ on small category $\BB$ a \em left Kan extension of $F$ along $Y$ \em is a functor $\mathit{Lan}_Y F \colon \BD \xra{ } \BC$ together with a natural transformation $\alpha \colon  F\Rightarrow  \mathit{Lan}_Y F \circ  Y$ such that for any other such pair $(G\colon \BD \xra{ } \BC,\beta \colon F\Rightarrow  G\circ Y)$ $\beta$  factors uniquely through $\alpha $. $\alpha$ is a natural isomorphism, if $Y$ is full and faithful; if $Y$ is a full embedding one even has $\alpha = id$. 
$\Lan_{ Y}F$ can be constructed with  $\Lan_{ Y}F$ acting  on a $\BD$-object  $D$ as   the colimit of the diagram $\mbox{$Y\/\/\/\downarrow\/\/D$}\xra{ \Phi_D}\BB\xra{F}\BC$ where $\Phi_D$ is the obvious forgetful functor into $\BB$, if  these colimits exist in $\BC$ for each $D$ (see e.g. \cite[Chapter X.3]{MacL}). 

In particular, if the  Yoneda embedding $Y_{\BB^\op}$ factors as $\BB^\op\xra{ Y}\BA\overset{I}\hookrightarrow\Set^\BB $ for some full subcategory $\BA$ of $\Set^\BB$ and   $\BC$ has $\mbox{$Y\/\/\/\downarrow\/\/A$}$-colimits for each $A$ in $\BA$,  one has natural isomorphisms 
\begin{equation}\label{eqn:Kanadjunct}
\BC(\mathit{Lan_{Y} } F(A),C)\simeq \Set^{\BB}(A, \hat{R}_FC)
\end{equation}
 where $\hat{R}_F \colon \BC\xra{ }\Set^{\BB}$ is the functor with $\hat{R}_F C  =\BB\xra{F^\op}\BC^\op\xra{ \BC^\op(C,-)=\BC(-,C)}\Set $ for each $\BC$-object $C$. Thus, the \em Yoneda extension \em $\Lan_{Y_{\BB^\op}}F$ has $\hat{R}_F$ as its right adjoint.

\item[\sf Sifted colimits.]\label{item4} A small category $\BD$ is called \em sifted \em if finite products in $\Set$ commute with colimits over $\BD$. Colimits of diagrams over sifted categories are called \em sifted colimits\em. Sifted colimits in a variety are precisely those colimits which are---like all limits---created by the forgetful functor. See \cite{ARV} for more details on this type of colimits.

  \item[\sf Equifiers.] Given an $I$-indexed family of pairs of functors $F^i,G^i\colon \BC\rightarrow \BC_i$ and for each $i\in I$ a pair of natural transformations $\phi^i,\psi^i\colon F^i\Rightarrow G^i$, we call (following \cite{AR}) the full subcategory of $\BC$ spanned by all objects $C$ satisfying $\phi^i_C = \psi^i_C$ for all $i\in I$ the \em joint equifier \em $\Equi(\phi^i,\psi^i)_i$ of the family $(\phi^i,\psi^i)$.
  
\item[\sf Some  notations.] The hom-functors of a category $\BC$ are denoted by $\BC(C,-)$. The equalizer of a parallel pair $(f,g)$ in $\BC$ is denoted by $\eq_\BC(f,g)$. If $(P,(P\xra{\pi_i}X_i)_i)$ is a product  and $(X\xra{ f_i}X_i)_i$ a family of morphisms in a category $\BC$ we denote by $X\xra{ \langle (f_i)_i\rangle}P$ the unique $\BC$-morphism with $\pi_i\circ  \langle (f_i)_i\rangle = f_i$ for each $i$. The \em $X$-fold copower \em of an object $C$ is denoted $X\cdot C$. Dually, given a coproduct $(C,(X_i\xra{{\mu_i}}C)_i)$, by  $C\xra{ [ (f_i)_i]} X$ we denote the unique morphism induced by a family of morphisms $(X_i\xra{ f_i}X)_i$. Functor categories are denoted by $\BB^\BA$ or $[\BA,\BB]$. If, for some category $\BC$, there is a category  $\BA(\BC)$ of ``structured $\BC$-objects" with forgetful functor $\BA(\BC)\xra{ |-|}\BC$ we will occasionally,  by slight abuse of notation,   denote an object in $\BA(\BC)$ with underlying $\BC$-object $C$ simply by $\mc$. 

We use the following notation for some frequently used categories:
 $\Set$ denotes the category of sets, $\Grp$ ($\Ab$), that of (abelian) groups,   $\Mon$ (${_c\Mon}$), that of (commutative) monoids, $\Ring$, that of unital rings, $_R\Mod$ ($\Mod_R$) that  of left (right) $R$-modules, $\mathit{Alg}_R$ ($_c\mathit{Alg}_R$), that of (commutative) $R$-algebras over a commutative ring. 
 
 For $n\in\N$ we often denote an $n$-element set simply by $n$.
\end{description}

\subsection{Varieties}
\subsubsection{Varieties by signatures}
A \em signature \em is a pair $\Sigma =(\bar{\Sigma}, \bar{\Sigma}\xra{ ar}\N)$ consisting  of a set $\bar{\Sigma}$ and an \em arity map \em $\bar{\Sigma}\xra{ ar}\N$. By $(\Sigma_n)_n$  we denote the family of preimages of $n$ under $ar$.  In the sequel we usually don't distinguish notationally between $\bar{\Sigma}$ and $\Sigma$.  $\Alg\Sigma$ then is the category of $\Sigma$-algebras  $\ma = (A,(A^{ar\sigma}\xra{ \sigma^\ma}A)_{\sigma\in\bar{\Sigma}})$ and $T(\Sigma)_n$ the set of $\Sigma$-terms in $n$ variables in the sense of universal algebra. $T(\Sigma)$, the disjoint union of all sets  $T(\Sigma)_n$,  is the set of  $\Sigma$-terms.

The interpretation of a term $t\in T(\Sigma)_n$ on a $\Sigma$-algebra $\ma$ will be denoted by $A^n\xra{ t^\ma}A$. 
Since the familiar definitions of $t^\ma$ are either given element-wise or use freeness of the $\Sigma$-algebras of $\Sigma$-terms and none of these methods works for arbitrary internal algebras in a category $\BC$ with finite products,  Freyd \cite{Freyd} suggested a recursive  definition using products, which however has some shortcomings as already observed in its review. 
We  therefore use the following recursive definition of  the  interpretation  $A^{n}\xra{ t^\ma}A$ of a $\Sigma$-term $t\in T(\Sigma)_n$ in a $\Sigma$-algebra $\ma$, which then will make sense  in any  category $\BC$ with finite products in Section \ref{sec:ia} below, where for the fundamental terms, that is the ``constants" $\sigma\in \Sigma_0$ and the ``variables" $x_1,\ldots,x_n,\ldots$ we add a subscript $n$ to indicate that they are considered as elements of $ T(\Sigma)_n$ (and not of some  $ T(\Sigma)_m$ with $m > n$) in order to avoid ambiguity: 
\begin{enumerate}
\item    $x_{n,j}^\ma:= 
A^n\xra{ {\pi_j}}A$ is the $j^{th}$ product projection, for the  variables $x_{1},\ldots,x_{n}$.
\item {$t^\ma :=  {A^n\xrightarrow{\langle (t_{i}^{\mathsf{A}})_i\rangle}A^{ar(\tau)}\xrightarrow{\tau^{\mathsf{A}}}}A$, if $t={\tau(t_1,\cdots,t_{ar(\tau)})}$ {with $t_i\in T(\Sigma)_n$}.  In particular 
\begin{enumerate}
\item $\sigma_n^\ma  := A^n\xra{!}1\xra{{\sigma^\ma}}A$ for each $t=\sigma \in \Sigma_0$,  where $!$ denotes  the unique morphism,
\item $\sigma(x_{1},\ldots,x_{n})^\ma  = A^n\xra{ \sigma^\ma}A$ for each $\sigma \in \Sigma_n$ and $n>0$.
\end{enumerate}
}
\end{enumerate}

Since every homomorphism $\ma\xra{ f}\ma'$ satisfies $t^{\ma'}\circ f^{n} = f\circ t^\ma$ for every $t\in T(\Sigma)_{n}$, every  such $t$ determines a natural transformation $\tilde{t}\colon |-|^n\Rightarrow |-|$ with components	 $\tilde{t}_\ma= t^\ma$, where $\Alg\Sigma\xra{ |-|}\Set$ denotes the forgetful functor.  By the Yoneda lemma the assignment $t\mapsto \tilde{t}$ defines a bijection between $T(\Sigma)_n$ and the set of natural transformations $nat(|-|^n,|-|)$, since $|-|^n$ is represented by the free $\Sigma$-algebra on $n$ whose carrier set is $T(\Sigma)_n$.  A \em set of $\Sigma$-equations \em is a family   $\BE = (\BE_n)_n$ with $\BE_n\subset T(\Sigma)_n\times T(\Sigma)_n$. As usual we denote  pairs  $(s,t)\in\BE$ by $s=t$. A $\Sigma$-algebra $\ma$ is said to \em satisfy  the equation $s=t$ \em iff $t^\ma = s^\ma$. $\Alg(\Sigma,\BE)$ denotes  the full subcategory of  $\Alg\Sigma$ of algebras satisfying all equations in $\BE$ and is called a \em variety\em. The pair $(\Sigma,\BE)$ is occasionally called an \em equational theory\em.

Writing $t = \tau(t_1,\ldots t_m)$ with $\tau\in\Sigma_m$ and $t_i\in T(\Sigma)_{k}$ and $s = \sigma(s_1,\ldots s_n)$ with $\sigma\in\Sigma_n$ and $s_j\in T(\Sigma)_{k}$   a $\Sigma$-algebra $\ma$ satisfies $s=t$ if and only if the following diagram commutes.
 \begin{equation}\label{diag:ax}
\begin{aligned}
\xymatrix@=2.5em{
    A^{k} \ar[r]^{ \langle (t^\ma_i)_i\rangle }\ar[d]_{ \langle (s^\ma_j)_j\rangle }&{A^m } \ar[d]^{ \tau^\ma}  \\
 A^n   \ar[r]_{ \sigma^\ma}          &    A 
}
\end{aligned}
\end{equation}
Satisfaction of an equational axiom $s=t$, thus, is the requirement {for  Diagram (\ref{diag:ax}) to commute}.
Since the paths of this diagram are the components of the natural transformations $\tilde{t}$ and $\tilde{s}$,  the variety  $\Alg(\Sigma,\BE)$ is,  as a subcategory of  $\Alg\Sigma$, the  joint equifier  $\Equi(\tilde{s},\tilde{t})_{s=t\in \BE}$.

\subsubsection{Varieties by Lawvere theories}\label{sec:Law}

The categorification of these notions is due  to Lawvere \cite{Law}, as is well known: Pairs $(\Sigma,\BE)$ are substituted by \em Lawvere theories \em $\T$ and varieties by the full subcategories $\Alg\T$ of $\Set^\T$ consisting of all finite product preserving functors.  
For the convenience of the reader not already familiar with the concept of Lawvere theory we add the following  remarks.
\begin{Definition}\label{def:lawv}\rm
A \em Lawvere theory \em is a category $\T$ with a countable set of objects $T_0, T_1,\ldots, T_n, \ldots$ (occasionally written simply as $0,1,\ldots, n,\ldots$) such that $\T$ has finite products and the object $T_1$ allows for each $n\in\N$  a specified product $(T_n,(T_n\xra{ \pi_{n,i}}T_1)_{1\leq i \leq n})$. Lawvere theories and finite product preserving functors constitute the category $\Law$.

The category $\Alg\T$ of \em $\T$-algebras \em has  its objects all finite product preserving functors $\T\xra{ A}\Set$ and its morphisms the natural transformations  between those. Evaluation at $T_1$ defines a faithful functor $\Alg\T\xra{ |-|}\Set$.
\end{Definition}

We recall the following basic facts of this categorification as follows, introducing some notations used in this note at the same time. For details we refer to \cite{ARV}.
\begin{enumerate}\renewcommand{\theenumi}{\Roman{enumi}}
\item  The dual of a skeleton of the category of finite sets is a Lawvere theory $\cn$. This is the initial object of $\Law$. For any Lawvere theory $\T$ we denote  by $T$ the unique $\Law$-morphism $\cn\xra{}\T$.
 There is  a forgetful functor from $\Law$ into the category of signatures   acting as $\T\xra{} \T(n,1) =:\Sigma_\T $. 
 This functor  has a left adjoint which assigns to a signature $\Sigma$  the \em free theory $\T_{\Sigma}$ over $\Sigma$\em. 
 In particular, $\Sigma_{\T_{\Sigma}} = T(\Sigma)$. There is a concrete isomorphism 
\begin{equation}\label{eqn:lift}
\Alg\T_\Sigma\simeq  \Alg\Sigma
\end{equation}
given on objects by $\mathsf{A}\mapsto (\mathsf{A}(T_1),(\mathsf{A}(\sigma(x_1,\cdots,x_{ar\sigma})))_{\sigma\in\Sigma})$.

\item There exists a dual biequivalence between Lawvere theories and varieties. In more detail: 
\begin{enumerate}
\item For every variety  $\BV=\Alg(\Sigma,\BE)$ there exists an essentially unique Lawvere theory $\TV$ such that $\BV$ is concretely equivalent to $\Alg\TV$. $\TV$ can be chosen to be the dual of the  full subcategory $\langle F_\BV n\mid n\in\N \rangle $  of $\BV$  formed by choosing, for each $n\in \N$, a free $\BV$-algebra $F_\BV n$  in such a way that $F_\BV n\neq F_\BV m$ for $n\neq m$. Equivalently, {$\TV(n,1)$ is the set of natural transformations $|-|^n\Rightarrow |-|$.}  In other words, $\T_\BV^\op$ can be considered as a full subcategory of $\Alg\T_\BV$.
\item For every Lawvere theory $\T$ the category $\Alg\T$ is concretely equivalent to some variety $\BV=\Alg(\Sigma,\BE)$, equivalently, $\T$ is isomorphic to $\TV$ for some $\BV$. 

\item Every morphism $\BS\xra{\Phi}\T$ of Lawvere theories, that is, a finite product preserving functor, induces the \em algebraic functor  \em      $\Alg\T \xra{\Phi^\ast}\Alg\BS$  given by $(\BT\xra{\ma}\Set)\mapsto (\BS\xra{\Phi}\BT\xra{ \ma}\Set)$.
This is a concrete functor admitting a left adjoint $\Phi_\ast$. Every concrete functor $\BV\xra{ U}\BW$ between varieties has a left adjoint $L_U$ and so determines a morphism in $\Law$:  Take the dual of the restriction $\langle F_\BW n\mid n\in\N \rangle \xra{L_U } \langle F_\BV n\mid n\in\N \rangle $.

In particular, for any Lawvere theory $\T$,  $\Alg\T\xra{T^\ast}\Alg\cn\simeq\Set$ is the forgetful functor;
it  creates limits and sifted colimits and reflects isomorphisms. Its left adjoint $T_\ast$ is the {\em free algebra functor}.
\end{enumerate}
\item Every Lawvere theory $\T$ is a regular quotient of a free one. In more detail: If $\Alg\T$ is concretely isomorphic to the variety $\Alg(\Sigma,\BE)$, then there exists  a quotient morphism $\T_\Sigma\xra{[-] }\T$ in $\Law$. This functor identifies terms $s$ and $t$, whenever $(s,t)$ is  deducible from $\BE$, that is, if  every algebra in  $\Alg(\Sigma,\BE)$ satisfies the equation $(s,t)$.  The functor $\Alg\T\xra{ [-]^\ast}\Alg\T_\Sigma$ is a full embedding.
In the language of universal algebra this quotient morphism corresponds to the canonical quotients $F_\Sigma n\xrightarrow{[-]}F_\BV n$ for each $n$ and, correspondingly, for each $\sigma\in \Sigma_n$, the $\T$-morphism $T_n\xra{[\sigma] }T_1$ can, in view of II.(a) above, be identified with the $\BV$-homomorphism $F_\BV 1\xra{ \bar{\sigma}}F_\BV n$ mapping 1 to $[\sigma(x_1,\ldots,x_n)]$.

\item There are essentially one-to-one correspondences $\T\mapsto \Bbb{T}$ between Lawvere theories and finitary monads on $\Set$, such that $\Alg\T\simeq\Set^\Bbb{T}$,  and the functors $\Phi^\ast $ induced by theory morphisms and  monad morphisms, respectively. The monad $\Bbb{T}$  is given by the adjunction $-\cdot F1\dashv \Alg\T(F 1,-)$.
\end{enumerate}

We will use below the following elementary construction of Lawvere theories. Let $A$ be an object in  a category $\BC$ with finite products. Choose, for each $n\in\N$, a power $A^n\xra{ \pi_i}A$ such that $A^n\neq A^m$ for $n\neq m$. Then the full subcategory of $\BC$ spanned by these powers is a Lawvere theory, which we will denote by $\T_\BC[A]$. Note that in this notation the theory $\T_\BV$ of a variety $\BV$ as described in II. (a) above is  $\T_{\BV^\op}[F_1]$.

We note moreover that, given an object $A$ in a category $\BC$ with coproducts, there is the monad $\Bbb{A}$ given by the adjunction 
$-\cdot A\dashv \BC(A,-)$. Though this monad formally is related to the  theory $T_{\BC^\op}[A]$ as $\Bbb{T}$ is related to 
$\T_{\BV^\op}[F_\BV1]$, $\Set^\Bbb{A}$ will not be concretely equivalent to a variety in general.

\subsubsection{Characterizing varieties as concrete categories}\label{sec:Linton} 
The characterization of varieties as concrete categories, that is, as categories $\BA$ equipped with faithful functor $\BA\xra{ U}\Set$ is essentially due to Linton \cite{Lin}. A concrete category $(\BV,U)$ is concretely equivalent to a variety if only if (0) {$\BV$} has coequalizers, (1) $U$ has a left adjoint, (2) $U$ preserves and reflects regular epimorphisms,  (3) $U$ preserves and reflects kernel pairs, and (4) $U$ preserves directed colimits. Assuming that $\BA$ even is cocomplete and has kernel pairs and regular factorizations of morphisms $U$ is a representable functor $U\simeq \BA(G,-)$ and one can replace the above conditions by the following ones:
$G$ is (a) a regular generator, (b) regular projective (i.e., projective with respect to regular epimorphisms), and (c) finitely presentable 
(see e.g. \cite[Thm 3.9.1]{Bor} or \cite{PorE}). We will call such $G$ a \em varietal generator\em.

One so obtains that for any  varietal generator $G$ in a cocomplete category $\BA$ with kernel pairs and regular factorizations of morphisms the variety $\BV$ with $\T_\BV =(\T_{\BA^\op}[G])^\op$ admits an equivalence $\BA\xra{ R_G}\BV$ with $|-|\circ R_G = \BA(G,-)$.

\subsubsection{A universal property  of varieties}\label{sec:sift}
Maybe somewhat lesser known than what we recalled in the previous section is the characterization of varieties $\Alg\T$ as free completion of the category $\T^\op$ under sifted colimits (see \cite{AR2}).  This is in detail: For  every functor $F\colon \T^\op\xra{ }\BC$  into a category $\BC$  with sifted colimits there exists an essentially unique sifted colimit preserving functor $S_F\colon\Alg\T\xra{ }\BC$   extending $F$ such that  the assignment $F\mapsto S_F$ defines an equivalence between  the categories $[\T^\op,\BC]$ and  $_{\text{\rm sift}}[\Alg\T,\BC]$, the category of sifted colimit preserving functors from $\Alg\T$ to $\BC$ (see \cite{AR2},\cite{ARV}). 

Since this result will be of importance  below we use the opportunity to describe the functor $S_F$ as the left Kan extension
$\Lan_Y(F)$ as follows, where $Y$ denotes the embedding $\T^\op\hookrightarrow\Alg\T$. 

\begin{Proposition}\label{prop:LR}
For every functor $F\colon \T^\op\xra{ }\BC$   into a category $\BC$  with sifted colimits  $\Lan_{Y}  F$ exists.
If $\BC$ in addition has finite coproducts, that is, if $\BC$ is cocomplete
the following hold.
\begin{enumerate}
\item $\Lan_{Y}  F$ coincides with the restriction of the Yoneda extension of $F$.
\item $\Lan_{Y}  F$ coincides with the unique sifted colimit preserving extension $S_F$ of $F$.
\item $S_F=\Lan_YF$  has a right adjoint $R_F$, if $F$ preserves finite coproducts.
\end{enumerate}
\end{Proposition}

The proof is based on the following lemma.
\begin{Lemma}[\cite{ARV}]
$\Alg\T$ is closed in $\Set^\T$ under sifted colimits. 
 For every $\T$-algebra $\T\xra{A}\Set$    the category $\El A$ of elements of the functor $\ma$ is sifted and coincides with  the comma category $\mbox{$Y\/\/\/\downarrow\/\/A$}$.
\end{Lemma}
\begin{Proof} (of Proposition \ref{prop:LR})
The first statement and item 1 are clear by the construction of Kan extensions.
By the lemma above the embedding $\Alg\T\xra{ I}\Set^\T$ preserves sifted colimits. Since the Yoneda extension is a left adjoint, $\Lan_YF$ preserves these by item 1 and, hence, coincides with $S_F$.
Concerning item 3 observe first, that the functor $\T\xra{ \hat{R}_F(C)}\Set$ preserves finite products if  $F^\op$ does; hence, in this case  $\hat{R}_F$  factors as  $\BC\xra{R_F }\Alg\T\xrightarrow{I}\Set^\T$.  By Equation \eqref{eqn:Kanadjunct} there are natural isomorphisms $\BC(\mathit{Lan_{Y}}  F(A),C) \simeq \Set^{\T}(I(A), \hat{R}_FC)\simeq\Alg\T(A, R_FC)$; in other words there is an adjunction $\Lan_YF\dashv R_F$.
\end{Proof}

The following diagram illustrates this situation. (Note, that the  Yoneda extension  $\Lan_{ Y_{\T^\op}}F$ will only exist if $\BC$ in addition has finite coproducts, that is, if $\BC$ is cocomplete.) We only include it in the diagram to indicate that in this case  the adjunction $\Lan_YF\dashv R_F$ appears as a restriction of the adjunction of the Yoneda extension).
$$
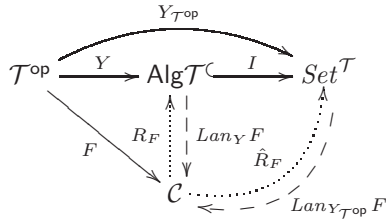
\xymatrix@=3em{ \T^\op \ar[dr]_{F}  
\ar@/^1.5pc/@{->}[rr]^{Y_{\T^\op} } \ar@{->}[r]^{Y} & \Alg\T \ar@{-->}@<1ex>[d]^{\Lan_Y F}
\ar@{^(->}[r]^{I} & \Set^\T\ar@/^1.8pc/@{-->}@<1ex>[dl]^{{\Lan_{ Y_{\T^\op}}F}}    \\
& \BC \ar@/_1.8pc/@{.>}@<0.5ex>[ur]^{\ \ \ \hat{R}_{{F}} } 
\ar@{.>}@<0.5ex>[u]^{R_F}   &  
}$$

\subsubsection{Internal algebras in a category}\label{sec:ia}

Given an arbitrary  category $\BC$ with finite products
 one defines categories $\Alg(\T,\BC)$, 
 $\Alg(\Sigma,\BC)$ and $\Alg((\Sigma,\BE),\BC)$    of \em internal algebras in $\BC$\em, equipped  with underlying functors $|-|_\BC$, by replacing  in the definitions of $\Alg\Sigma$, $\Alg(\Sigma,\BE)$ and $\Alg\T$, respectively,  ``set" by ``$\BC$-object" and ``map" by ``$\BC$-morphism". 
   For a variety $\BV = \Alg(\Sigma,\BE)$ one calls $\Alg((\Sigma,\BE),\BC)$ the category of \em internal $\BV$-algebras in $\BC$\em. 
  A $\BV$-algebra in $\BC$, thus, is a pair $(A,(\sigma^A)_{\sigma\in\Sigma})$ such that for each $(s,t)\in \BE$ Diagram \eqref{diag:ax}, considered as a diagram in $\BC$, commutes.
 Internal $\BV$-algebras in $\BC$ can alternatively be described as follows (\cite[Chapter III.6]{MacL},\cite{Freyd}):
 
 \begin{Fact}\label{fact:intalg}\rm
 A $\BC$-object $A$ is (the underlying object) of an internal $\BV$-algebra iff, for each $\BC$-object $C$, the set $\BC(C,A)$ is (the underlying set) of a $\BV$-algebra.

This is equivalent to saying that $\BV$-algebras $\ma$ in $\BC$ with underlying object $A$ correspond essentially one-to-one to functors $R_\ma$ making the following diagram commute.
\begin{equation}\label{Txlift}
\begin{aligned}
\xymatrix@=2.5em{
{ \BC}\ar[r]^{R_\ma } \ar[dr]_{\BC(-,A)}&{\BV} \ar@{->}[d]^{|-| }     \\
 &{ \Set}     
}
\end{aligned}
\end{equation}
\end{Fact}

 The following facts either follow by direct generalization or are obvious. 
\begin{Fact}\label{prop:fund}\rm
If $\T_\BV$ is the Lawvere theory of a  variety  $\BV=\Alg(\Sigma,\BE)$, then the following hold.
\begin{enumerate}
\item\label{item:1!} The functor $\Alg(\T_\BV,\BC)\xra{ |-|_\BC}\BC$ creates limits and reflects isomorphisms. \\
For every $t\in\T_\Sigma(n,1)$ there is a natural transformation $\tilde{t}\colon |-|_\BC^n\Rightarrow |-|_\BC$, whose components are the  $\BC$-morphisms $ \ma(t)\colon|\ma|^n_\BC\xra{}|\ma|_\BC$. 
\item\label{item:ce} $\Alg(\TV,\BC)$ is concretely equivalent over $\BC$ to $\Alg((\Sigma,\BE),\BC)$, the  full subcategory $\Equi(\tilde{s},\tilde{t})_{s=t\in \BE}$ of the category $\Alg(\Sigma,\BC)$. 
\item  Morphisms $\BS\xra{\Phi}\T$ of Lawvere theories determine functors  $_\BC\Phi^\ast\colon \Alg(\T,\BC) \lra\Alg(\BS,\BC)$   given by $(\BT\xra{\ma}\BC)\mapsto (\BS\xra{\Phi}\BT\xra{ \ma}\BC)$.
\item\label{factS} Every product preserving functor $\BC\xra{ S}\BD$ between categories with finite products maps $\T$-algebras to $\T$-algebras, that is, $S$ can be lifted to a  functor $\Alg(\T,\BC)\xra{ ^\T\!{S}}\Alg(\T,\BD)$ by mapping a $\T$-algebra $\ma$ to $S\circ\ma$. ${ ^\T\! S}$ preserves all limits which are preserved by $S$.
  In the language of universal algebra, one has a functor $\Alg(\BV,\BC)\xrightarrow{^\BV\! S}\Alg(\BV,\BD)$ such that $^\BV\! S(\mathsf{A})=(S(A),(S(\sigma^{\mathsf{A}}))_{\sigma\in \Sigma})$.
\end{enumerate}
However, {unlike} to the case of $\BC=\Set$,  neither the forgetful functor $\Alg(\T,\BC)\xra{|-|_\BC }\BC$ nor the functors $_\BC\Phi^\ast$ have left adjoints in general.  They do so, if $\BC$ is a locally finitely presentable category, and in this case $\Alg(\T,\BC)$ is locally finitely presentable (see e.g. \cite{Por}).
\end{Fact}

\section{Internal coalgebras}\label{sec:IC}
\subsection{Descriptions of the category of coalgebras}\label{sec:def}

If $\BC$ is a category with finite coproducts one calls an internal algebra in $\BC^\op$ an \em internal coalgebra in $\BC$\em. This generalizes the terminology used already   by Kan \cite{Kan} as early as 1958 in the special situation, where $\BC=\Mon$  and $\BV=\Grp$.  In this section we describe various equivalent descriptions of this concept in some detail.
  
\subsubsection{Coalgebras for a variety}\label{sec:eq}
 \begin{Definition}\label{def:var}\rm
Let  $\BC$ be a  category with finite coproducts.

  For any Lawvere theory $\T$ the category of  \em internal $\T$-coalgebras in $\BC$   \em  is the category \linebreak$\Coalg(\T,\BC) =\Alg(\T,\BC^\op)^\op$. 

For any variety  $\BV =\Alg(\Sigma,\BE)$ the category   {\em internal $\BV$-coalgebras in $\BC$} (or \em internal $(\Sigma,\BE)$-coalgebras\em)  is the category $\Coalg(\BV,\BC) =\Alg(\BV,\BC^\op)^\op$.
\end{Definition}

 \begin{Facts}\label{fact:equif}\rm
 \begin{enumerate}
\item $\Coalg(\T,\BC)$ is the category of finite coproduct preserving functors $\T^\op\xra{ }\BC$. 
\item $\Coalg(\BV,\BC) = \Coalg((\Sigma,\BE),\BC)$ is,  by dualization of Facts \ref{prop:fund}.\ref{item:ce},  the joint equifier  $\Equi(\hat{s},\hat{t})_{s=t\in \BE}$ in $\Coalg(\Sigma,\BC)$,  where  $\hat{t}:=\tilde{t}^\op\colon ||-||\Rightarrow n\cdot ||-||$ is the natural transformation  
determined by $t\in\T_\Sigma(n,1)$.  $t_\ma:= \hat{t}_\ma\colon A\to n\cdot A$ is the interpretation of $t$ in an internal  $\Sigma$-algebra in $\BC^\op$, for each such $t$.  Moreover, $\Coalg(\BV,\BC) \simeq \Coalg(\T_\BV,\BC)$.
\end{enumerate}
\end{Facts}
 
 Adopting the language of universal algebra, a $\Sigma$-coalgebra in $\BC$ is a  pair $\ma = (A, (A\xra{ \sigma_\ma}n\cdot A)_{\sigma\in\Sigma_n})$, where $A $ is a $\BC$-object and $A\xra{ \sigma_\ma}n\cdot A$ is a $\BC$-morphism, called the \em $n$-ary co-operation determined by $\sigma\in \Sigma_n$\em. A $(\Sigma,\BE)$-coalgebra in $\BC$ then is a $\Sigma$-coalgebra in $\BC$  
making all  Diagrams (\ref{diag:ax})  in $\BC^\op$ commute. This is often expressed by saying that    it \em  satisfies duals of the equational axioms \em defining  $\BV$, which  means explicitly that each of the following diagrams commutes in $\BC$.
 \begin{equation*}
\begin{aligned}
\xymatrix@=2.5em{
 A\ar[r]^{ \sigma_\ma}\ar[d]_{ \tau_\ma}  & {n\cdot A }\ar[d]^{ [{((s_j)}_\ma)_j] }    \\
 m\cdot A \ar[r]_{ [{((t_i)}_\ma)_i] }          &    k\cdot A 
 }
 \end{aligned}
\end{equation*}
The following diagrams from \cite{Kan}, with $C\xra{ \mu}C+C$ the comultiplication, $C\xra{ \epsilon}0$ the counit, and $C\xra{\iota}$ the coinversion,  displaying the axioms of a cogroup may serve as an illustration.
\begin{equation*}
\xymatrix@=2.5em{ C\ar[r]^\mu\ar[d]_\mu&{C+ C}\ar[d]_{[\mu,\id_C]}\ar@/^1pc/@{.>}[d]^{\mu+ \id_C}\\
C+ C\ar[r]^{\!\!\!\!\!\!\!  [\id_C,\mu]} \ar@/_1pc/@{.>}[r]_{ \id_C+\mu}&C+ C+ C}
\ \ \ \ \ \ \ \ \ \ \ \   
\xymatrix@=2.5em{ 0+C\ar[d]_{\nu_2^{-1}}&C\ar[d]_\mu
{\ar[d]_\mu}\ar[r]^{\nu_1}\ar[l]_{\nu_2}&C+0\ar[d]^{\nu_1^{-1}}\\
C &C+ C\ar[r]_{\ \  [\id_C,\epsilon]}\ar[l]^{\!\!\!\!\!\!\!\  [\epsilon,\id_C]}&C}
\end{equation*}
\begin{equation*}
 \begin{minipage}{5cm}
 \xymatrix@=2.5em{
 C+C \ar@{.>}[d]_{\nabla =[\id_C,\id_C]}& C+C\ar[dl]^{[\id_C,\iota]} \ar@{.>}[l]_{\id_C +\iota}&
 C\ar[d]_\epsilon \ar[r]^{\!\!\!\!\! \mu} \ar[l]_{\ \  \mu}& C+C \ar@{.>}[r]^{\iota+ \id_C}\ar[dr]_{[\iota,\id_C]}   & C+C\ar@{.>}[d]^{[\id_C,\id_C]=\nabla}\\
C && 0\ar[rr]_{!} \ar[ll]^{!} & &C
}
\end{minipage}
 \end{equation*}


\begin{Remark}\label{rem:}\rm
It may happen that $\Coalg(\BV,\BC)$ is the category $\mathsf{1}$. This is the case for example if $\BC(C,0) =\emptyset$ for all $C\neq 0$ (as in $\Set$) and $\BV$ has no constants. An example is $\Coalg(\Mon,\Set)$.

It also may happen that $\Coalg(\BV,\BC)$ is isomorphic to the category $\BC$. Examples are the categories  $\Coalg(\Mon,\BC)$ and $\Coalg(_c\Mon,\BC)$, if $\BC$ is a semi-additive category.  This is due  to the familiar fact that the (only) co\-monoids  in a cartesian monoidal category $\BC$ are the triples $(C,\Delta_C,!_C)$ with $C\xrightarrow{\Delta_C}C\times C$  the diagonal and $C\xrightarrow{!_C}1$  the unique morphism, and that these are cocommutative.

Other examples are $\Coalg(\Grp,\BC)$ and  $\Coalg(\Ab,\BC)$ if $\BC$ is additive. In this case, for each $\BC$-object $C$, the  hom-functor $\BC(C,-)$ can be lifted to a functor $\BC\to\Ab$. In other words, $C$ carries the structure of an $\Ab$-coalgebra in $\BC$. By the above the claim follows and one obtains equivalences of categories $\Coalg(\Grp,\BC)\simeq \BC\simeq \Coalg(\Ab,\BC)$.
\end{Remark}

\subsubsection{Coalgebras as representable functors}\label{sec:repfct}

\begin{Definition}\label{def:repres}\rm
Given a  Lawvere theory $\T$ and  $\BC$ a  category with finite coproducts. 
 A functor $\BC\xra{R_A }\Alg\T$ is called \em $\T$-representable (on $\BC$) by $A\in \ob\BC$\em\ or an \em $\Alg\T$-lift of $\BC(A,-)$\em, if  $\BC\xra{ }\Alg\T\xra{ |-|}\Set = \BC \xra{\BC(A,-)}\Set$ for some $\BC$-object $A$, that is, if diagram
 \eqref{Txlift} commutes. 
 
 $\mathsf{Rep}(\T,\BC)$ denotes the category of $\T$-representable functors  on $\BC$ with all  natural transformations as its morphisms.

Occasionally we will rather talk about  $\BV$-lifts and $\BV$-representable functors, if $\T =\T_\BV$ for a variety  $\BV$.  
\end{Definition}
\begin{Remarks}\label{rem:repconv}\rm
\begin{enumerate}
\item  For any morphism $R_A\overset{\mu}\Rightarrow R_{A'}$ in  $\mathsf{Rep}(\T,\BC)$ the natural transformation $\BC(A,-)\overset{|\mu|}\Rightarrow \BC(A',-)$ determines by the Yoneda Lemma a $\BC$-morphism $A'\xra{ f_\mu}A$ between the representing objects. $ \mathsf{Rep}(\T,\BC)^\op$ then is concrete over $\BC$ by means of   the {faithful} functor 
  $\{-\}_\T$  given by $(R\xra{ \mu} R' )\mapsto (A'\xra{ f_\mu} A)$.
 \item For any $\BV$-representable functor 
  $R_A$ one has for each set $X$  the natural isomorphism
\begin{equation}\label{eqn:natrep}
\BV(FX,R_A-)\simeq\BV(F1,R_A-)^X\simeq \BC(A,-)^X\simeq\BC(X\cdot A,-).
\end{equation}
Identifying $\BV(FX, R_A-)$ and $\BC(X\cdot A,-)$, for each $\BV$-homomorphism $F1\xra{ g}FX$ the natural transformation $\BV(g,R_A-)\colon \BV(FX,R_A-)\Rightarrow \BV(F1,R_A-)$ is a natural transformation \linebreak $\BC(X\cdot A,-)\Rightarrow \BC( A,-)$ and so determines by the Yoneda lemma a $\BC$-morphism $X\cdot A \xra{\gamma}A$. 
\end{enumerate}
\end{Remarks}
 
Dualizing Fact \ref{fact:intalg} one obtains

\begin{Lemma}\label{lem:1}
For every Lawvere theory $\T$ and every category $\BC$ with finite coproducts the assignment $\ma\mapsto R_\ma$ defines an  
equivalence, concrete over $\BC$,  $$\Coalg(\T,\BC)\simeq\mathsf{Rep}(\T,\BC)^\op.$$
\end{Lemma}

In more detail: 
If $\T^\op\xra{\ma}\BC$ is a coalgebra with $\ma(T_1) = A$, then the functor $\BC(A,-)$ has the $\T$-lift $\BC\xra{ R_\ma}\Alg\T$ given by the assignment $C\mapsto \T\xra{\ma^\op }\BC^\op\xra{ \BC(-,C)}\Set$. 
Conversely, for any $\T$-representable functor $R_A$ one obtains  a coalgebra $\T^\op\xra{ \ma }\BC$ with $\ma T_1= A$ by the
 assignment $F1\xra{ g}Fn\mapsto A\xra{ \gamma}n\cdot A$ defined in item 2. above.
For a more conceptual description of this correspondence see below.

 \begin{Remark}\label{rem:11}\rm
A representable functor $\BC(A,-)$ may have non-isomorphic $\Alg\T$-lifts $H,H'\colon\BC\xra{ }\Alg\T$, equivalently, a $\BC$-object may carry more than one structure of a $\T$-coalgebra.

For example  there are the functors $\Phi, \Psi\colon \mathit{Ring}\lra\Mon$ mapping a unital ring  $R$  to its additive monoid $(R,+,0)$ and its multiplicative monoid $(R,\times,1)$, respectively. These are $\Mon$-lifts of the representable functor $\Ring(\Z[x],-)$, which is up to equivalence the forgetful functor of $\Ring$. See Theorem \ref{corr:algfu} for a complete discussion of this phenomenon.
\end{Remark}

\subsubsection{Coalgebras as adjunctions}\label{sec:adj}

For every adjunction $L\dashv R\colon \Alg\T\xra{ }\BC$, where $\BC$ is a  category with finite coproducts, the right adjoint $R$ is $\T$-representable. In fact, the natural isomorphism $\BC(LF1,-)\simeq \Alg\T(F1,R-)$ is equivalent to $|-|\circ R\simeq\BC(LF1,-)$. The converse holds, if $\BC$ is cocomplete. This  has been claimed already, for the case of $\BC$ being cocomplete,  by Freyd (without giving a proof) and is
 an immediate consequence of Beck's lifting theorem for adjunctions along a monadic functor (see e.g. \cite{Bor}). 
For the sake of the reader interested in an algebraically flavoured  construction of this left adjoint 
we explain below how this argument works in the given  situation, correcting at the same time the proof given in \cite{Berg2}. 

The resulting correspondences between coalgebras, representable functors and adjunctions can be seen more conceptually as follows, where  $\Adj(\BD,\BC)$ denotes the category of adjunctions from $\BD$ to $\BC$ as defined in \cite{MacL}: 

\begin{Theorem}\label{theo:pres}
For every Lawvere theory $\T$ and every cocomplete category $\BC$  there are  equivalences 
 $$\Coalg(\T,\BC)\simeq \mathsf{Rep}(\T,\BC)^\op\simeq  \Adj(\Alg\T,\BC).$$
 In more detail:
 \begin{enumerate}
\item  The left adjoint corresponding to a $\T$-coalgebra $\ma$ in $\BC$ is its left Kan extension $\Lan_Y\ma$. It is the sifted colimit preserving extension of $\ma$ as well.
\item The $\T$-representable functor corresponding to a $T$-coalgebra $\ma$ in $\BC$ is the right adjoint of  $\Lan_Y\ma$  
\item  The $T$-coalgebra $\ma$ in $\BC$ corresponding to a left adjoint functor $\Alg\T\xra{ L}\BC$ is its  restriction to $\T^\op$, considered as a subcategory of $\Alg\T$. 
\end{enumerate}
\end{Theorem}
Indeed, by Proposition \ref{prop:LR} there is the adjunction $L_\ma= Lan_Y\ma\dashv R_\ma$ for every $\T$-coalgebra $\ma$ in $\BC$.  By its definition the functor $R_\ma$ is the $\T$-representable functor represented by $\ma(T_1)$ \tr({see Section \ref{sec:repfct}).  
 Noting further that the equation $S_F=\Lan_YF$ implies that the correspondence of item 3 is an equivalence by the definition of the free cocompletion (see Section \ref{sec:sift})}, this theorem only   rephrases Proposition \ref{prop:LR}.

\begin{Remarks}\label{rems:f}\rm
\begin{enumerate}
\item $\mathsf{Rep}(\T,\BC)$ is  equivalent  to the category $\Radj(\BC,\Alg\T)$ of right adjoints from $\BC$ to $\Alg\T$ with natural transformations as morphisms, while 
$\Adj(\Alg\T,\BC)$ is  equivalent  to the category $\Ladj(\Alg\T,\BC)$ 
 of left adjoints (equivalently, by the Special Adjoint Functor Theorem, to ${_{coc}[\Alg\T,\BC]}$, the category of cocontinous functors) from $\Alg\T$ to $\BC$ with natural transformations as morphisms. 
 $\Adj(\Alg\T,\BC)$ becomes concrete over $\BC$  by assigning to  a morphism $L\dashv R\xra{\langle \sigma,\rho \rangle} L'\dashv R'$, that is,  a  pair of  conjugate  transformations, the $\BC$-morphism  $L(F1)\xra{\sigma_{F1} }L'(F1)$,  the $\sigma$-component at the free $\T$-algebra over a singleton. 
\item It is easy to see that for any natural transformation $L \overset{\sigma}{\Rightarrow}L'$ the underlying $\BC$-morphism of the right adjoint $R$ of $L$ is $\{\rho\}_\T$, where $\langle \sigma,\rho \rangle$ is a conjugate pair. Thus, the equivalences of the theorem are concrete over~$\BC$.
\item Denoting by $\BV$ the variety corresponding to $\T$, the $\BV$-coalgebra $\ma$ in $\BC$ corresponding to an adjunction $L\dashv R\colon \BV\to \BC$ then is, in the language of universal algebra, $\mathsf{A}=(L(F_\BV 1),(\sigma_{\mathsf{A}})_{\sigma\in\Sigma})$, with co-operations $\sigma_{\mathsf{A}}=L(F_\BV 1)\xrightarrow{L(\bar{\sigma})}L(F_\BV n)$, $\sigma\in \Sigma_n$ (see Item III. of Section  \ref{sec:Law}). 
\end{enumerate}
\end{Remarks}

\subsubsection*{\sf The Eilenberg-Watts Theorem.}
The theorem above can be seen as the utmost generalization of the Eilenberg-Watts Theorem, since combining it with Facts  \ref{fact:equif} we 
 the get the equivalence
$$\Coalg(_R\Mod,{_S\Mod})\simeq  \Ladj(_R\Mod,{_S\Mod)}.$$ 
It thus suffices to show that $\Coalg(_R\Mod,{_S\Mod})$ is nothing but the category ${_S\Mod_R}$  of left $S$-, right $R$-bimodules. We prove something more in order to show how the respective result from \cite{N/S} fits into this picture. Let $\BC$ be a cocomplete and additive category.
 Then, for every $\BC$-object $C$ the set $\BC(C,C)$ carries the structure of a ring $\End_\BC(C)$ and the categories $\Coalg(\Ab,\BC)$ and $\Ab$ are isomorphic  (see Section \ref{sec:eq});  hence,  $\Coalg(_R\Mod,\BC)$ is isomorphic to the category of \em left $R$-objects in $\BC$ \em in the sense of \cite{Pop}, that is, of pairs $(C,\rho)$, where $R\xra{ \rho}\End_\BC(C)$ is a ring homomorphism. 
In particular the category $\Coalg(_R\Mod,{_S\Mod})$ is nothing but the category ${_S\Mod_R}$  of left $S$-, right $R$-bimodules. Note that this argument strengthens the original result in that the assumption of additivity is not needed.
(In fact, it wasn't needed from the very beginning since every finite coproduct preserving functor between additive categories is additive.).

 \subsubsection*{\sf An algebraic proof of right adjointness of $\BV$-representable functors.}
Choose,  for any algebra $A$ in a variety $\BV$, 
its \em canonical presentation by generators and relations \em given by
the coequalizer diagram $\xymatrix@=2em{
    F|K_A| \ar@<.5ex>[r]^{ r_A }\ar@<-.5ex>[r]_{ s_A }&{F|A| } \ar[r]^{\epsilon_A } &A }$, 
where $\epsilon_A$ is  the counit  of the adjunction $F\dashv |-|$, $\xymatrix@=2em{
    K_A \ar@<.5ex>[r]^{ p }\ar@<-.5ex>[r]_{ q }&{F|A| } \ar[r]^{\epsilon_A } &A 
}
$ 
its congruence relation (kernel pair) and $(r_A,s_A)$ is the pair of homomorphic extensions of the projections $p$ and $q$. Denote by $\rho_A$ and $\sigma_A$ the $\BC$-morphisms $|A|\cdot R\xra{ }|K_A|\cdot R$ corresponding to $r_A$ and $s_A$ 
 according to Remarks \ref{rem:repconv}.2. and let $\xymatrix@=2em{
    |A| \cdot R\ar@<.5ex>[r]^{ \rho_A }\ar@<-.5ex>[r]_{ \sigma_A }&{|K_A|\cdot R } \ar[r]^{\ \ \ \upsilon_A } &LA 
}$ be a  coequalizer diagram.
Since contravariant hom-functors map coequalizers to equalizers one obtains, for each  $C$ in $\BC$,
the following equalizer diagrams  in $\Set$.\vspace{-2ex}
\begin{eqnarray*}\small
\xymatrix@=3em{
{\BV(A,HC) }\ar[r]^{\!\!\!\!\!\!\BV(\epsilon_A,HC) } & \BC(|A|\cdot R,C)\ar@<.5ex>[r]^{\BC(\rho_A,C)}\ar@<-.5ex>[r]_{\BC(\sigma_A,C)}& \BC(|K_A|\cdot R,C)
}\\
\xymatrix@=3em{
{\BC(LA,C) }\ar[r]^{\!\!\!\!\!\!\BC(\upsilon_A,C) } & \BC(|A|\cdot R,C)\ar@<.5ex>[r]^{\BC(\rho_A,C)}\ar@<-.5ex>[r]_{\BC(\sigma_A,C)}& \BC(|K_A|\cdot R,C)
}
\end{eqnarray*}\vspace{-1ex}
Hence, for each $C$ in $\BC$ and each $A$ in $\BV$, there is a bijection  $\BC(LA,C)\xra{\alpha_{A,C}}\BV(A,HC)$
\begin{equation}\label{eqn:crux}\small
 \begin{aligned}
\xymatrix@=3em{
{ \BC(LA,C)}\ar@{.>}[d]_{\alpha_{A,C}}\ar[r]^{ \!\!\!\!\!\!\BC(\upsilon_A,C)} & \BC(|F|A||\cdot R,C)\\
{ \BV(A,HC)} \ar[ur]_{\ \ \ \BV(\epsilon_A,HC) }          &
}
\end{aligned}
\end{equation}
The claim now follows by \cite[Chapter IV.1, Corollary 2]{MacL}, 
provided that these bijections are natural in $C$. But this is obvious: Since $\BV(\epsilon_A,H-)$ and 
$\BC(\upsilon_A, -)$ are natural and $\BV(\epsilon_A,HC)$ is monic, for each $C$ in $\BC$, commutativity of Diagram \eqref{eqn:crux} implies that the family $(\alpha_{A,C})_C$ is natural. It is this naturality condition of which  in the proof in \cite{Berg2} is not (and cannot be) taken care of,  due to the fact that there is only chosen an arbitrary presentation by generators and relations of the algebra $A$ at the beginning. 

\subsubsection{Coalgebras for monad}\label{sec:monad}
Recall that a $\BC$-object $A$ carries the structure of a $\T$-coalgebra in $\BC$ iff there exists a morphism of theories $\T\xra{ \phi}\T_{\BC^\op}[A]$. If $\BC$ is cocomplete, such $\phi$ extends to     a monad morphism from the free algebra monad $\mathbb{T}$ of $\Alg\T$ to the monad $\mathbb{A}$ given by the adjunction $-\cdot A\dashv \BC(A,-)$. 

Since monad morphisms $\mathbb{T}\xra{\phi }\mathbb{A}$ correspond bijectively to functors $\Set^\mathbb{A}\xra{ \phi^\ast}\Set^\mathbb{T}$ commuting with the forgetful functors (see e.g. \cite[Section 3.6]{BW}) we obtain a bijective correspondence between such pairs $\ma=(A,\phi)$ and \em $\mathbb{T}$-representable functors \em $\BC\xra{ R_\ma}\Set^\mathbb{T}$, that is, functors $R_\ma$ such that $U_\mathbb{T}\circ R_\ma = \BC(A,C)$, for each $C$ in $\BC$, and these are nothing but $\T$-representable functors, since  $\Set^\mathbb{T}\simeq\Alg\T$.

Hence,  a $\BC$-object $A$ carries the structure of a $\T$-coalgebra $\ma$ in $\BC$ iff there exists a morphism of monads $\Bbb{T}\xra{ \phi}\Bbb{A}$, and the $\T$-representable functor $R_\ma$ corresponding to $\ma$ is nothing but $\BC\xra{ K_A}\Set^\Bbb{A}\xra{ \phi^\ast}\Set^\Bbb{T}\simeq\Alg\T$ where $K_A$  is the comparison functor of the adjunction $-\cdot A\dashv \BC(A,-)$.
Thus, the following definition specializes for finitary monads to the case discussed so far and allows for an obvious generalization of Theorem \ref{theo:pres} from finitary monads on $\Set$ to arbitrary ones.
\begin{Definition}\label{def:}\rm
Let $\BC$ be a category with arbitrary coproducts and $\mathbb{T}$ be a monad on $\Set$. A \em $\mathbb{T}$-coalgebra in $\BC$ \em is a pair $\ma=(A,\phi)$, where $A$ is $\BC$-object  and $\mathbb{T}\xra{ \phi}\mathbb{A}$ a monad morphism.\\
A morphism $(A,\phi)\xra{f} (B,\psi)$ of $\mathbb{T}$-coalgebras is a $\BC$-morphism $A\xra{ f}B$ such that $\BC(f,-)$ lifts   to a natural transformation $R_\mb\Rightarrow R_\ma$, that is, if  for each $C$ the map $\BC(B,C)\xra{ \BC(f,C)}\BC(A,C)$ is  a $\mathbb{T}$-morphism.
\end{Definition}

It is easy to see that some of the theoretic results below hold, if the assumption \em $\T$ is a Lawvere theory \em would be replaced by \em $\Bbb{T}$ is a monad on $\Set$ \em and, correspondingly $\Alg\T$ by $\Set^\Bbb{T}$. This applies in particular for Proposition \ref{prop:props}, Theorem \ref{corr:algfu}, and Proposition \ref{theo:adj}. 

\subsection{Some examples, counter examples and applications}\label{sec:ex}

\begin{description}
\item[\sf Cogroups.] 
For historical reasons only we start with the probably first example of cogroups in a  category with finite coproducts $\BC$, that is of $\Grp$-coalgebras in $\BC$, from homotopy theory (see \cite{Kan},\cite{Freyd}). Let $\BC = \mathit{HTop}$ be the category of pointed topological space as objects and homotopy classes $\overline{f}$ of base-point preserving continuous maps $f$ as morphisms. The coproduct of two pointed space $(X,x_0)$ and $(Y,y_0)$ in $\mathit{HTop}$ is their wedge sum $X\vee Y$, the topological sum of $X$ and $Y$ with $x_0$ and $y_0$ identified and taken as its base point. 

By definition the $n^{th}$ homotopy group $\pi_n(X,x_0)$  of a pointed space $(X,x_0)$ has as its underlying set the hom-set $\mathit{HTop}((S^n,s_0),(X,x_0))$. Thus, the functor $\mathit{HTop}\xra{ \pi_n}\Grp$ is a $\Grp$-representable functor with representing object~$S^n$. The corresponding cogroup structure is described as follows:

When collapsing an equator through $s_0$ of a pointed $n$-sphere $(S^n,s_0)$ to $s_0$, the quotient is  $S^n\xra{ \mu}S^n\vee S^n$. The multiplication of $\pi_n(X,x_0)$ is given by $\overline{f}\star\overline{g} = \overline{[f,g]\circ\mu}$, the unit is $1 = \overline{{x_0}}$ with ${x_0}$ the constant map with value $x_0$, and the inversion is  given by $\overline{f}\mapsto \overline{f\circ \iota}$ where $\iota$ is a change of orientation of $S^n$. Then proving that these data define a group is equivalent to proving that $(S^n, s_0)$ is a cogroup  in $\mathit{HTop}$ with comultiplication  $\mu$ and inversion $\iota$.

\item[\sf Modules.]\label{ex:can03} Let $\BV$ be the variety $\Mod_R$ for a commutative ring and $M$ an $R$-module. 
The internal hom-functor $[M,-]\colon \Mod_R\xra{ }\Mod_R$ is $\BV$-representable with representing object $M$ and has a left adjoint $M\otimes -$. The corresponding coalgebra has as its underlying object the module $M$ and the co-operations $M\xra{+}M+ M $ and $M\xra{r_{M}}M$ for each $r\in R$, which are by Theorem \ref{theo:pres} the homomorphisms
$M\simeq M\otimes R\xra{M\otimes + } M\otimes R^2 \simeq M^2$ and $M\simeq M\otimes R\xra{M\otimes r\cdot- } M\otimes R \simeq  M$, respectively. Thus, the coalgebra structure on $M$ is trivial in that it only ``adds'' the diagonal.

\item[\sf Group algebras.]
Consider the varieties $\BV= {\mathit{Alg}_R}, \ {_c\BV}={_c\mathit{Alg}_R}, \ \BW = \Grp, \ {_c\BW} =\Ab$,  $\mathcal{M} =\Mon, \ {_c\mathcal{M}} = {_c\Mon} $ and ${\mathcal{U}}=\Mod_R$.  
The following diagram displays the (obvious) forgetful functors, where $\Phi^\ast(A)$ is the additive group of an algebra $A$, while $\Psi^\ast(A)$ is its  multiplicative monoid, and the functor $(-)^\times$  assigns  to a (commutative) monoid its (abelian) group of invertible elements.
\begin{equation}\label{diag:12}
\begin{aligned}
\xymatrix@=2em{
 \Mod_R   \ar[d]_{\Xi^\ast}  &  \mathit{Alg}  \ar@{->}[dl]_{\Phi^\ast}  \ar[r]^{ \Psi^\ast } \ar[l]_{ \Sigma^\ast  }& \Mon\ar[r]^{(-)^\times}  & \Grp\ar@/_2pc/@{->}[ll]^{R[-]} \\
\Ab          &    {_c\mathit{Alg}_R}\ar[r]^{_c\Psi^\ast }
\ar@{^(->}[u]
& {_c\Mon}\ar@{^(->}[u]\ar[r]^{_c(-)^\times} & \Ab\ar@{^(->}[u]
}
\end{aligned}
\end{equation}
Being algebraic all forgetful functors have a left adjoint;  the functors $(-)^\times$ and ${_c}(-)^\times$ are right adjoints of the respective embeddings $\Lambda^\ast$ and  ${_c\Lambda}^\ast$. The left adjoint of $(-)^\times\circ \Psi^\ast$ is the group algebra functor $R[-]$ (note that $R[-] = \Psi_\ast\circ E$ where $E$ is the embedding $\Grp\hookrightarrow \Mon$). The adjunction $R[-]\dashv (-)^\times\circ \Psi^\ast\colon \Grp\xra{ } \mathit{Alg}_R$ hence determines a cogroup in $\mathit{Alg}_R$. In the sequel we will omit the subscript $_c$ at functors if no confusion is possible.

\item[\sf Hopf structures.]\label{Hopf} 
Since the monoidal structure of ${_c\mathit{Alg}_R}$ coincides with the cocartesian structure, that is, $\otimes = +$ in $_c\mathit{Alg}_R$,  the  comonoids in $_c\mathit{Alg}_R$ are the commutative  bialgebras over $R$, that is, $\Coalg(\Mon, {_c\mathit{Alg}_R})\simeq {_{c}\Bialg_R}$, and $\Coalg(_c\Mon, {_c\mathit{Alg}_R})\simeq {_{bi}\Bialg_R}$, the category of bi-commutative bialgebras.
Similarly, the  cogroups in $_c\mathit{Alg}_R$ are the commutative  Hopf algebras over $R$, that is, $\Coalg(\Grp, {_c\mathit{Alg}_R}) = {_{c}{\Hopf_{\!R}}}$, while $\Coalg(\Ab, {_c\mathit{Alg}_R})\simeq{_{bi}{\Hopf_{\!R}}}$, the category of bi-commutative Hopf algebras over $R$.

Denoting by   $\Ab\xra{ \Lambda^\ast}{_c\Mon}$ the forgetful functor, the embedding ${_{bi}{\Hopf_{\!R}}}\hookrightarrow  {_{bi}\Bialg_R}$ is given by the functor ${_{bi}{\Hopf_{\!R}}}\simeq \Coalg(\Ab, {_c\mathit{Alg}_R})\xra{^{{_c\Alg_R}\!}\Lambda} \Coalg(_c\Mon, {_c\mathit{Alg}_R})\simeq {_{bi}\Bialg_R}$ and, thus,   has a right adjoint by Proposition  \ref{prop:props} below. This improves a result of \cite{POR}, where this fact has been shown for the case of absolutely flat rings only. 

Obviously  the adjunction $R[-]\dashv (-)^\times\circ \Psi^\ast$ of the previous example restricts to an adjunction 
$^\Ab{R[-]}\dashv {(-)}^\times\circ {\Psi}^\ast\colon \Ab\xra{ }{_c\mathit{Alg}_R}$  and so determines a bicommutative Hopf algebra. The familiar Hopf algebra structure on a commutative group algebra $R[G]$ then is given by evaluation of the functor 
$\Coalg(\Ab,\Ab)\xra{ {^{\Ab}\!R[-]}} \Coalg(\Ab,{_c\mathit{Alg}_R}) \simeq {_{bi}\Hopf}$  at the abelian group $G$, since
$\Coalg(\Ab,\Ab)\simeq\Ab$ by the Eckmann-Hilton argument. 
\item[\sf A counter example.]
The full subcategory $\mathit{Tor}$ of $\Ab$ spanned by all commutative torsion groups  is a coreflective subcategory of (hence a comonadic category over)  $\Ab$ but it is not of the form $\Coalg(\T,\Ab)$.

Assume   $\mathit{Tor} =\Coalg(\T,\Ab)$ for some Lawvere theory $\T$. Representing $\T$  as
a regular quotient of some $\T_\Sigma$, we can consider   $\Coalg(\T,\Ab)$ as an equifier in the functor category $\Coalg Q_\Sigma$ for the functor $Q_\Sigma X = \prod_{n\in\N}  (n\cdot X)^{|\Sigma_n|} = \prod_{n\in\N}   X^{|n\Sigma_n|}$
according to Fact \ref{fact:equif}. Since $\Ab$ has biproducts the functor  $\Ab\xra{ Q_\Sigma}\Ab$ preserves  products, 
such that the forgetful functor $\Coalg Q_\Sigma\xra{ }\Ab$ creates  products and the equifier corresponding to $\Coalg(\T,\Ab)$ is closed under  products. But  $\mathit{Tor}$ fails to be   closed under products in $\Ab$. 
\end{description}

\subsection{Some properties of $\Coalg(\T,\BC)$}\label{sec:prop}
\begin{Proposition}\label{prop:props}
The following hold for any  cocomplete category  $\BC$  and any Lawvere theory  $\T$. 
\begin{enumerate}
\item $\Coalg(\T,\BC)$ has all colimits which exist in $\BC$ and these are created by $|\!|-|\!|$.
\item Every hom-functor $\Coalg(\T,\BC)\xra{ }\Set$ has a left adjoint.
\item For every $\BS\xra{\Phi}\T$ the functor 
 $\,{^\BC}\Phi:= \Coalg(\T,\BC)\xra{_{\BC}{\Phi^\ast}^\op}\Coalg(\BS,\BC)$ has a right adjoint~${_\BC\Phi}$. 
 \item The functor  $\Coalg(\T,\BC)\xra{{ ^\T\!{S}}}\Coalg(\T,\BD)$ with respect to a finite coproduct preserving functor $\BC\xra{ S}\BD$ defined as in Facts\ref{prop:fund}.\ref{factS}   has a  right  adjoint ${_\T S}$, provided that $S$ preserves colimits. 
\end{enumerate}
\end{Proposition} 
\begin{Proof}
Item 1 follows by dualization of Fact \ref{prop:fund}.\ref{item:1!} and, hence Item 2 follows  trivially. 
The right adjoint of Items 3 can  be constructed by (the dual of) the familiar lifting theorem of adjunctions (see e.g \cite[Chapter 4.5]{Bor}). In particular, the unit of the adjunction ${^\T\!S}\dashv {_\T{S}}$ (occasionally we denote this adjunction by ${^\BV\!S}\dashv {_\BV{S}}$ if $\T$ is the theory of the variety $\BV$) is point wise an equalizer. 
The same holds for the adjunction ${^{\BC}\Phi}\dashv {\Phi_{\BC}}$.
\end{Proof}

Bergman\cite{Berg1}  studied further categorical properties of $\Coalg(\T,\BV)$ for varieties $\BV$ using highly technical arguments; he  showed in particular that $\Coalg(\T,\BV)$ is a complete category and its underlying functor into $\BV$ has a right adjoint. In  this section we improve these results by the use of standard arguments from the theory of locally presentable categories. Hence,  for the rest of this section the base category $\BC$ in which $\T$-coalgebras are formed is assumed to be to be a locally presentable category. 
\begin{Theorem}\label{thm:crux} 
For any  Lawvere theory  $\T$  and any locally presentable category $\BC$  the category  $\Coalg(\T,\BC)$ is locally presentable.
Its forgetful functor $\Coalg(\T,\BC)\xra{|\!|-|\!|}\BC$ has a right adjoint $C_{\T\!,\BC}$ and, moreover, is comonadic.
\end{Theorem}

\begin{Proof}
The claim follows from the equivalence $\Coalg(\T,\BC)\simeq  \Ladj(\Alg\T,\BC)$, since for locally presentable categories  $\BC$ and $\BD$ the category $\Ladj(\BC,\BD)$ is locally presentable (see Section \ref{not}). 

The existence  of a right adjoint now follows from Item 1 of the preceding proposition  by  (the dual of) the Special Adjoint Functor Theorem because their domains are locally presentable. Comonadicity follows by the Beck-Par\'e-Theorem: the category  $\Coalg(\T,\BC)$, considered as a subcategory of the (comonadic) category $\Coalg Q_\Sigma$ (see the counter example of Section \ref{sec:ex} for the definition of $Q_\Sigma$), is closed under absolute equalizers. 
\end{Proof}

\begin{Example}\label{ex:}\rm
The cofree functor $C_{\T\!,\BC}$ can occasionally be constructed. As an example  consider the case where $\BC = {_S\Mod}$ and $\T$ is the theory of ${_R\Mod}$, hence, $\Coalg(\T,\BC)$ is the category $_S\Mod_R$ of bimodules. Since this is a variety $\BV$, the functor $C_{\T\!,\BC}$ needs to be a $\BV$-representable functor, whose representing $S$-module  has, somehow, to encode the ring $R$. A natural choice, thus, is  the $S$-module $F|R|$, where ${_R\Mod}\xra{ |-|}\Ab$ denotes the forgetful functor and $\Ab\xra{ F}{_S\Mod}$ the ``free'' functor, that is, the left adjoint of the forgetful functor $\Mod\xra{ }\Ab$, by abuse of notation also denoted by $|-|$. By the natural equivalence ${_S\Mod}(F|R|,N)\simeq\Ab(|R|,|N|)$ it suffices to show that $\Ab(|R|,|N|)$ belongs to   $_S\Mod_R$, for every $S$-module $N$. But this obvious with  $(s\cdot f)(r):=s\cdot (f(r))$ and $(f\cdot r)(r'):=f(rr')$. 

We note that, since $_S\Mod_R$ is also equivalent to $\Alg(\Mod_R,{_S\Mod})$, by Fact~\ref{prop:fund} the functor  $|\!|-|\!|$ also has a left adjoint. As easily seen this is given by $M\mapsto |M|\ot_\Z R$, where the abelian group $|M|\ot_\Z R$ 
inherits a left $S$-action from that of $M$, given by $s\cdot (y\otimes z):=(s\cdot y)\otimes z$, and  a right $R$-action  from multiplication in $R$, $(x\otimes y)\cdot r:=x\otimes (yr)$. 
\end{Example}

Kan already had constructed comonoids in groups on each free group \cite{Kan}, while Freyd has shown that every free algebra in a variety $\BV$ carries the structure of a $\BW$-algebra, provided there exists an algebraic functor $\BV\xra{ }\BW$ \cite{Freyd}. We now describe this construction in greater generality as follows, where $\BC\xra{ F}\Alg(\T,\BC)$ is the free algebra functor which exists by Fact \ref{prop:fund}. 

\begin{Proposition}\label{prop:can}
Let $\BC$ be a locally presentable category  and $\T$ a Lawvere theory. Then there exists a functor $\BC\xra{T}\Coalg(\T,\Alg(\T,\BC))$ such that  $\BC\xra{T}\Coalg(\T,\Alg(\T,\BC))\xra{|\!|-|\!| }\Alg(\T,\BC) = \BC\xra{ F}\Alg(\T,\BC)$. $T$ has a right adjoint.
\end{Proposition}
\begin{Proof}
By the dual of Lemma \ref{lem:1} the functor $\BC(-,A)$ factors as $\BC\xra{ } \Alg\T\xra{|-| }\Set$ iff $A=|\ma|$ for a (essentially) unique $\T$-algebra $\ma$ or, in other words, for each $\BC$-object $C$ the set $\BC(C,A)$ is the underlying set of a (essentially) unique $\T$-algebra $\ma_C$. From the natural bijection $\Alg(\T,\BC)(FC,\ma)\simeq \BC(C,|\ma|_\BC)$  one now concludes that for each $\T$-algebra $\ma$ in $\BC$ the set $\Alg(\T,\BC)(FC,\ma)$ is the underlying set of the  $\T$-algebra $\ma_C$, such that the assignment $\ma \mapsto \ma_C$ defines a functor $\Alg(\T,C)\xra{R_C }\Alg\T$ which is a $\T$-lift of  $\Alg(\T,\BC)(FC,-)$. Define $T$ by $C\mapsto R_C$.

Since $F$ preserves and $|\!|-|\!|$ creates colimits, the functor $T$ preserves them. Thus, the final statement follows by the Special Adjoint Functor Theorem whose assumptions on the category $\BC$ are satisfied by every locally presentable category.
\end{Proof}

\section {Coalgebras in varieties}\label{sec:CV}

As from now we assume the base category $\BC$ in which $\T$-coalgebras are formed to be a variety. $\T$-coalgebras in $\Alg\BS$ are occasionally  called \em$\T\!$-$\,\BS$-bimodels \em or \em$\T\!$-$\,\BS$-bialgebras\em\ with the category of those denoted by $[\T,\BS]$. ($\T\!$-$\,\T$-bialgebras will in the sequel simply called \em $\T$-bialgebras\em.) 
 
 These categories can be used to define a bicategory having as objects all Lawvere theories,  the categories $[\T,\BS]$ as categories of 1-cells   with horizontal composition  defined by $(\BS\overset{\mb}{\Rightarrow}\BR)\odot (\T\overset{\ma}{\Rightarrow}\BS):= \BT^\op\xra{\ma}\Alg\BS\xra{L_\mb}\Alg\BR$ and  the embeddings $1_\T := Y^{\T}\colon\T^\op\hookrightarrow \Alg\T$ as units, 
 and  with natural transformations as 2-cells (see \cite[Chapter 15]{ARV}). Note that the ``product" $\mb\odot\ma$ can be seen, equivalently,  as the $\T\!$-$\,\BR$-bialgebra representing the composition $R_\ma\circ R_\mb$ of the representable functors $R_\ma$ and $R_\mb$, represented by $\ma$ and $\mb$, respectively. This construction is a straightforward  generalization of the familiar bicategory $\mathit{RING}$,  having as objects all unital rings (equivalently, all theories of module categories),  as categories of 1-cells $[R,S]$ the categories ${_R\Mod_S}$ of $R$-$S$-bimodules with $\ot$ as horizontal composition, and as 2-cells the bimodule homomorphisms.

We prefer to work with  the following 2-category $\mathit{LAWV}$, which is in view of Theorem \ref{theo:pres} biequivalent to the bicategory  just described, and where the categories of 1-cells can, alternatively, be chosen  as $\mathsf{Rep}(\Alg\T,\Alg\BS)^\op$.
\begin{enumerate}
\item objects are all Lawvere theories,
\item  the categories of 1-cells are the categories $\Ladj(\Alg\T,\Alg\BS)$ with composition of functors as horizontal composition and identities as units, 
\item and 2-cells are natural transformations.
\end{enumerate}

\subsection{Coalgebras in arbitrary varieties}\label{sec:co-var}
\subsubsection{Canonical constructions} 

As follows from  Proposition \ref{prop:can} every free algebra $F_\BV X$ in a variety $\BV$ carries the structure of a $\BV$-coalgebra. Obviously, one then obtains a $\BW$-coalgebra on $F_\BV X$ as well, if there is an algebraic functor $\BV\xra{ \Phi^\ast}\BW$.
We now will describe these structures explicitly.

Since every algebraic functor $\Phi^\ast$  induced by a morphism of Lawvere theories $\TW\xra{ \Phi}\TV$ 
commutes with the underlying functors, it is $\BW$-representable with representing object $F_\BV1$, the free $\BV$-algebra on one generator. We denote this $\BW$-coalgebra in $\BV$ by $V_\Phi(1)$. With  notations as in item 3 of Section \ref{sec:ex} and item 2 of Section \ref{sec:prop} one has $V_\Phi(1)= {^{\BV}\Phi}(V_\id(1))$. Considered as an  adjunction $V_\Phi(1)$ is nothing but $\Phi_\ast\dashv {\Phi^\ast}$. 
These are the only $\BW$-coalgebras on $F_\BV 1$. Since $\Coalg(\BW,\BV)$ has coproducts we can form the coalgebra $V_\Phi(X) = X\cdot V_\Phi(1)$ for each set $X$. Considered as an  adjunction this is $X\cdot\Phi_\ast\dashv {\Phi^\ast}^X$ with the $\BW$-representable functor  ${\Phi^\ast}^X$, represented by $F_\BV X$. This construction is  functorial  since every map $X\xra{ f}Y$ determines a natural transformation $\Phi^{\ast^Y}=\BV(F_\BV Y,-)\Rightarrow \BV(F_\BV X,-)=\Phi^{\ast^X}$, which is by the Yoneda lemma,  that $F_\BV(f)$ is a coalgebra morphism $V_\Phi(X)\xra{ }V_\Phi(Y)$. We so have got the following result, where the first statement is already contained in  \cite{Freyd}, and the last one follows trivially from the definition of $V_\Phi$ as the $X$-fold copower of $V_\phi(1)$ in $\Coalg(\BW,\BV)$.

\begin{Theorem}\label{corr:algfu}
There is an essentially one-to-one correspondence between
 morphism of Lawvere theories $\TW\xra{ \Phi}\TV$  and  $\TW$-coalgebra  $\BV_\Phi(1)$ on $F_\BV 1$;  moreover
\begin{enumerate}
\item the assignment $X\mapsto X\cdot V_\Phi(1)$ defines a functor $V_\Phi\colon\Set\xra{ }\Coalg(\BW,\BV)$ with $|\!|-|\!|\circ V_\Phi =F_\BV$,
\item $V_\Phi$ has the functor  
 $G_\Phi: = \hom(V_\Phi(1),-)$  as its right adjoint.
\end{enumerate}
\end{Theorem}

\begin{Remarks}\label{rem:19}\rm
\begin{enumerate}
\item With $V:=V_{\id}$ the following equivalences then are obvious.
\begin{equation}\label{commdiag}
{^{\BV\!}\Phi}\circ V \simeq V_\Phi\simeq {^\BW\!\Phi_\ast}\circ W
\end{equation}
\item The  co-operations of $\BV_\Phi( X)$ are the $\BV$-homomorphisms $F_\BV X\xra{ X\cdot \Phi(\sigma)}X\cdot F_\BV n \simeq n\cdot F_\BV X$ for all $\sigma\in \T(n,1) =\BW(F_\BW1, F_\BW n)$ and $n\in\N$. These co-operations  act on the free generators $x\in X$ of $F_\BV X$ as $X\xra{ \eta_X} F_\BV X\xra{\langle \nu_k\rangle} (F_\BV (n\cdot X))^n\xra{\Phi(\sigma)^{F_\BV (n\cdot X)}}  F_\BV (n\cdot X)\simeq n\cdot F_\BV X$, where $\eta_X$  is the insertion of generators and $\langle\nu_k\rangle$ is the homomorphism whose coordinates are the coproduct injections $F_\BV X\xra{ \nu_k}n\cdot F_\BV X$. In other words,
\begin{equation}\label{eqn:coopgen}
{\sigma}_{V_\Phi(X)}   \circ  \eta_X = \Phi(\sigma)^{F_\BV (n\cdot X)}\circ\langle \nu_k\rangle \circ  \eta_X
\end{equation}
\end{enumerate}
\end{Remarks}

This situation is illustrated by the following diagram,  where the unlabelled arrows are the forgetful functors. Here all cells obviously commute, while  commutativity  of the outer frame represents Equation \eqref{commdiag}.
\begin{equation*}\label{diag:Diag}
\begin{aligned}
\xymatrix@=1.2em{
 \Coalg(\TV,\BV)    \ar[rr]^{^{\BV\!}\Phi}\ar[dr]_{  }&&{\Coalg(\TW,\BV) }\ar[dl]   && \Coalg(\TW,\BW)\ar[ll]_{{^\BW\!\Phi_\ast}}\ar[dl]^{ }\\
&\BV            & &   \BW\ar[ll]_{ \Phi_\ast}\\
&& \Set \ar[ul]^{F_\BV}\ar@/^2pc/@{->}[uull]^{V} \ar@/_2pc/@{->}[uurr]_{W}\ar[ur]_{F_\BW}
}
\end{aligned}
\end{equation*}

A concrete description of the right adjoint of  $V_\Phi$ is obtained as follows where, for $\tau\in \T_\BW(n,1)$,  we denote by $\tau_\ma$  the underlying map of the $\BV$-homomorphism $\ma\tau$, by $\Phi(\tau)^{n\cdot A}$  the interpretation  of  $\Phi\tau$ in the $\BV$-algebra $n\cdot A$, and by   $\nu_1,\ldots, \nu_n\colon A\xra{ }n\cdot A$ in $\BV$ the  coproduct injections. 
\begin{Proposition}\label{theo:adj}
For every $\T_\BW$-coalgebra $\ma$ in $\BV$ the set $G_\Phi(\ma)$ is the subset 
\begin{equation}\label{eq:G}
G_\Phi \ma \  = \bigcap_{t\in\T_\BW(n,1),n\in\N} \eq_\Set(t_\ma, (\Phi(t))^{n\cdot A}\circ \langle \nu_k\rangle)
\end{equation}
of $| |\!|\ma |\!| |$ and the embeddings $e_\ma$ form are a natural transformation $e\colon G_\Phi\Rightarrow |-|\circ |\!|-|\!|$.
\end{Proposition} 

\begin{Proof} 
Identifying $\Coalg(\T_\BW,\BV)$ and $\Adj(\Alg\T_\BW,\BV)$,    the hom-set 
$\hom(V_\Phi(1),\ma)$ is ${nat}(\Phi_\ast,L_\ma)$, the set of natural transformations  from $\Phi_\ast$ into the 
 left adjoint $L_\ma$ of $\ma$ corresponding to $\ma$, for each $\T_\BW$-coalgebra $\ma$.
Given a natural transformation $\mu\colon\Phi_\ast\Rightarrow L_\ma$, one has $|\!|\mu|\!| = \mu_1\colon \Phi_\ast(F_\BW 1)= F_\BV 1\xra{ }L_\ma F_\BW 1 = \ma (1) =|\ma| $; hence, the assignment $\mu\mapsto \mu_1(1)$ defines an injective map $G_\Phi(\ma)\xra{ e_\ma}| |\!|\ma |\!| |$ and these maps form a natural transformation, trivially.

Equation \eqref{eq:G}  is equivalent to the following statement: For every  $\BV$-homomorphisms $F_\BV 1\xra{ f} A=\ma(1) $ the family $(n\cdot f)_n$ is  a natural transformation $\Phi_\ast\Rightarrow L_\ma$ if and only if  $f(1)\in G_\Phi \ma$. But this  is evident by the following diagram, using Equation \eqref{eqn:coopgen}. 
  \vspace{-0.5em}
\begin{equation*}\label{diag:}
\begin{minipage}{5cm}
\xymatrix@=2em{
   1  \ar[r]^{ \eta_1 }\ar[d]_{ \eta_1 }&{F 1} \ar[r]^{t_{V_\Phi(1)}}\ar[d]_{\langle \iota_k\rangle } & n\cdot F1\ar[d]^{id}\ar@/^2pc/@{->}[dd]^{n\cdot f} \\
 F1 \ar[r]^{\langle \iota_k\rangle}\ar[d]_f&  (Fn)^n  \ar[d]_{(n\cdot f)^n} \ar[r]^{\Phi(t )^{Fn}}          &  Fn \ar[d]^{n\cdot f}\\
A\ar@/_2pc/@{->}[rr]_{t_\ma}\ar[r]_{\!\!\!\!  \langle \nu_k\rangle} & (n\cdot A)^n\ar[r]_{\ \ \Phi(t )^{n\cdot A}}   & n\cdot A
}
\end{minipage}\vspace{-5ex}
\end{equation*}
\end{Proof}

\begin{Remarks}\label{rem:GS}\rm
The following facts are easy to verify.
\begin{enumerate}
\item  $G_\Phi \ma$ contains all $\BV$-constants in $A$.
\item By Equation  \eqref{eqn:coopgen}  for every set $X$ the following holds:  $X\subset G_\Phi V_\Phi(X)\subset |F_\BV X|$. In other words  the set $G_\Phi V_\Phi(X)$ generates the $\BV$-algebra $|\!|V_\Phi(X) |\!|$. In general this is not the case for arbitrary coalgebras and this implies that $G_\Phi$ in general fails to be faithful.
\item If $\T_\Sigma\xra{ [-]}\T_\BW$ is a regular quotient, then   
$$G_\Phi \ma = \bigcap_{\sigma\in\Sigma_n,n\in\N} \mathrm{Eq}_\Set(\sigma_\ma, (\Phi\sigma)^{n\cdot A}\circ \langle \nu_k\rangle).  \vspace{-0.5em}$$
\item Writing $G_\BV$ instead of $G_{\id_\BV}$  commutativity of Diagram \eqref{diag:Diag} implies by composition of adjunctions
\begin{enumerate}
\item The forgetful functor $\BW\xra{ |-|}\Set$ factors as $\BW\xra{ C_{\T_\BW,\BW}}\Coalg(\T_\BW,\BW)\xra{ G_\BW}\Set$ and 
\item  $\Coalg(\T_\BW,\BV)\xra{G_\Phi }\Set$ factors as $\Coalg(\T_\BW,\BV)\xra{ _{\T_\BW}{\Phi_\ast}}\Coalg(\T_\BW,\BW)\xra{ G_\BW}\Set$ and as $\Coalg(\T_\BW,\BV)\xra{^{\BV\!}\Phi}\Coalg(\T_\BV,\BV)\xra{ G_\BV}\Set$.
\end{enumerate}
\end{enumerate}
\end{Remarks}

\begin{Lemma}\label{lem:comp}
Given a  factorization $\Phi= \T_\BW\xra{\Xi }\T_\BU \xra{ \Sigma} \TV$   in $\Law$,  the following diagrams commute.
\begin{equation*}\label{eqn:}
 \begin{aligned}
\xymatrix@=1.5em{
{\Coalg(T_\BW,\BU) }\ar[rr]^{^\BW\Sigma_\ast } && \Coalg(\T_\BW,\BV)\\
& { \Set} \ar[ul]^{ U_\Xi} \ar[ur]_{V_\Phi }         
}
\end{aligned} \ \ \hfill\ \ 
\begin{aligned}
\xymatrix@=1.5em{
{\Coalg(T_\BW,\BU) }\ar[dr]_{ G_\Xi}   && \Coalg(\T_\BW,\BV)\ar[ll]_{_\BW\Sigma_\ast }\ar[dl]^{G_\Phi } \\
& { \Set}        
}
\end{aligned}
\end{equation*}
\end{Lemma}

\begin{Example}\label{grp-l}
With notations as in Diagram \eqref{diag:12} $ \Coalg(\Mon,{}_c\mathit{Alg}_R)$ is equivalent to the category ${}_c\mathit{Bialg}_R$  of all commutative unital and counital $R$-bialgebras, while $ \Coalg({\Ab},{}_c\mathit{Alg}_R)$ is equivalent to the category ${}_{bi}\Hopf_{\!R}$ of all bicommutative $R$-Hopf algebras.
\\
By Equation \eqref{eq:G}  one gets, for any  commutative bialgebra $\mathsf{A}=(A,m,e,\Delta,\epsilon)$, 
\begin{equation*}\label{PE}
G_{\Phi}(\mathsf{A})=\{\, a\in A\mid \Delta(a)=a\otimes 1 + 1\ot a\} 
\end{equation*}
that is, ${G}_{\Phi}(\mathsf{A})$ is the usual set of  primitive  elements of $\mathsf{A}$, while for any   bicommutative bialgebra 
\begin{equation*}\label{GLE}
 G_{\Psi}(\mathsf{A})=\{\, a\in A\mid \Delta(a)=a\otimes a,\ \epsilon(a)=e(1)\,\} 
\end{equation*}
is the usual set of  group-like elements of  $\ma$.

For ${_R\Mod_R}\simeq \Coalg(_R\Mod,{_R\Mod})$, the category of $R$-bimodules, and any bimodule $M$  
$$G(M) = \{\, m\in M\mid rm=mr  \ \  \forall  r\in R \}$$
is the space of $R$-invariants of $M$ in the sense of \cite{Ag}.
\end{Example} 

\begin{Remark}[Comonoids in the category of groups]\label{rems:Kan}\rm
With $\BV = \Grp$,  $\BW = \Mon$ and  the  forgetful functor $\Phi^\ast\colon \BV\xra{}\BW$   the following has been shown by Kan (see \cite{Kan}):   Every comonoid  in $\Grp$ has as its underlying group a free group $FX$ and every free group carries precisely one comonoid structure (\cite[Thm. 3.10]{Kan}). This theorem states in our terminology (with $X_\ma:=G_\Phi \ma\setminus \{1\}$): {\em  For every comonoid $\ma$ in $\Grp$ there is an isomorphism } $\ma\simeq V_\Phi(X_\ma)$.
 
$\ma$  then is an internal cogroup in $\Grp$ (\cite[Ex. 6.5]{Kan}), which is in our terminology $V(X_\ma)$; 
this is clear by Remarks \ref{rem:19}: the co-inversion of any cogroup $V(X)$ is   the homomorphic extension $FX\xra{\iota}FX$ of the map $X\xra{ }FX$ with $x\mapsto x^{-1}$, where $x^{-1}$ denotes the inverse of $x$ in the free group $FX$. 

This result implies  that the embedding $\Phi^\ast$  is the only algebraic functor from $\Grp$ to $\Mon$ and  the identity is the only algebraic functor on $\Grp$.  Combining Kan's result with the methods developed above we can give 
concrete descriptions of the following constructions of which we so far only know (see Theorem \ref{thm:crux}) that they exist (cp. \cite[pp 399-400]{Berg1}):
\begin{enumerate}
\item  A right adjoint  of the forgetful functor $\Coalg(\Grp,\Grp)\xra{|\!|-|\!| }\Grp$ is given by the functor $\Grp\xra{ |-|}\Set\xra{ V}\Coalg(\Grp,\Grp)$. 
In fact, the counit of the free group adjunction $|\!|V(|A|)|\!| =F|A|\xra{ \epsilon_A}A$ is the required counit: For any group homomorphism  $|\!|\mb |\!| =|\!|V(X_\mb) |\!|=  F (X_\mb)\xra{ f}A $ the homomorphism $F X_\mb\xra{ F(f)}F |A| $ is the only homomorphism $f^\ast$ with $\epsilon_{|A|}\circ f^\ast = f$ and coincides with $|\!|V(|f|) |\!|$. 
\item Products in  $\Coalg(\Grp,\Grp)$ can be constructed as follows,   since $\Coalg(\Grp,\Grp)\xra{ G}\Set $ is a right adjoint of $V$ and, hence, preserves products (while $V$ preserves coproducts):

Let $\ma$ and $\mb$ be cogroups in $\Grp$. Writing $\ma=V(X_\ma)$ and $\mb=V(X_\mb)$ we have 
 $X_{\ma\times\mb} =  G_\BV(\ma\times \mb) -1=G_\BV\ma\times G_\BV\mb -1= 
 (X_\ma+1)\times (X_\mb+1)-1$ and,  thus, $$\ma\times\mb =V(X_\ma\times X_\mb + X_\ma +X_\mb )=V(X_\ma\times X_\mb)+\ma+\mb.$$
\end{enumerate}
Kan's result is a rather special property of the variety $\Grp$ and its embedding into $\Mon$. 
As follows  from Remark \ref{rem:11} and  the discussion of linear groups in Section \ref{sec:3.2} below, in general   there exist coalgebras on non-free algebras as well, and  free algebras may allow for more than  one coalgebra structure.
The examples mentioned in Example \ref{ex:example1} below share with $\Grp$ the property that the coalgebras $V(X)$ in $\BV$ are the only $\BV$-coalgebras in $\BV$; as the argument above shows  also in these cases $|\!|-|\!| $ is left adjoint to $V\circ |-|$. In view of Item 2 of Remarks \ref{rem:GS} they provide however the other extreme case of the possible sizes of the sets $G_\BV V(X)$: while for $\BV=\Grp$ the set $G_\BV V(X)$ is the smallest possible one, in those cases it is the largest.
\end{Remark}

\subsubsection{Some applications}\label{sec:3.2}

\begin{description}
\item[\sf Change of rings.] 
Consider   a ring homomorphism $R\xra{ \phi}S$ (equivalently, a morphism of Lawvere theories $\Phi\colon \T_{_R\Mod}\xra{ }\T_{_S\Mod}$). 

The algebraic functor $\Phi^\ast\colon\BV= {_S\Mod}\xra{ }{_R\Mod}$ is the ${_R\Mod}$-represen\-table functor ${_S\!\hom(S,-)}$, known as the  \em restriction of  scalars-\em functor.

The corresponding $S$-left, $R$-right bimodule structure on $S$ is $V_\Phi(1)$. By the above its left adjoint $\Phi_\ast$ is the functor $V_\Phi(1)\ot_{R^{}} -$, known as the \em extension of scalars-\em functor.

Since the forgetful functors from module categories into $\Ab$ create colimits, the functor $\Phi^\ast$, being a concrete functor over $\Ab$,  preserves them. Consequently, the restriction of $\Phi^\ast$ to the dual of the theory of ${_S\Mod}$ is an ${_S\Mod}$-coalgebra in ${_R\Mod}$ and $\Phi^\ast$ is  its left Kan extension extension along the embedding ${\T_{_S\Mod}}^\op\hookrightarrow {_S\Mod}$,  hence a left adjoint,  its right adjoint $G$ being the ${_S\Mod}$-lift of $_R\hom(\Phi^\ast S,-)$. In other words, $G$ is the \em coextension of scalars-\em functor.

We so have obtained, without any calculations, the familiar adjunctions \\
\mbox{}\hfill  $\mathit{extension\, of \, scalars} \dashv  \mathit{restriction\, of \, scalars}\dashv \mathit{coextension\, of \, scalars}.$\hfill\mbox{}
\end{description}

\begin{description}
\item[\sf Morita theory.]
Lawvere theories $\T$ and $\BS$ are called \em Morita equivalent \em if  the categories $\Alg\T$ and $\Alg\BS$ are equivalent (not necessarily concretely so). Since an equivalence $L$ is a left and  a right adjoint of  its equivalence inverse $R$ 
one concludes from Theorem \ref{theo:pres} that the following are equivalent, where the respective coalgebras $\ma$  and $\mb$ are related by the property  $L_\ma \simeq R_\mb$    and  $L_\mb\simeq R_\ma $.
\begin{enumerate}
\renewcommand{\theenumi}{\roman{enumi}}
\item $\Alg\T\xra{ L}\Alg\BS$ is an equivalence with  equivalence inverse $\Alg\BS\xra{ R}\Alg\T$.
\item There exists a $\T$-coalgebra $\ma$ in $\Alg\BS$ with $L\simeq L_\ma$ and $R\simeq R_\ma$.
\item There exists a $\BS$-coalgebra $\mb$ in $\Alg\T$ with $R\simeq L_\mb$ and $L\simeq R_\mb$.
\end{enumerate}
Thus, Lawvere theories $\T$ and $\BS$ are  Morita equivalent, if they are equivalent in the
2-category $\mathit{LAWV}$.

The coalgebras $\ma$ and $\mb$ above can in view of  Section \ref{sec:Linton} equivalently be described as follows: $\ma (1)$ is a varietal generator in $\Alg\BS$ such that  $\T\simeq (\T_{{(\Alg\BS)}^\op}[\ma (1)])$ and $\mb (1)$ is a varietal generator in $\Alg\T$ such that  $\BS\simeq (\T_{{(\Alg\T)}^\op}[\mb (1)])$. This follows from the facts that  $\ma (1)= L(F_{\Alg\T}1)$ is a varietal generator in $\Alg\BS$ since $F_{\Alg\T}1$ is a varietal generator in $\Alg\T$ and equivalences clearly preserve varietal generators, and $\ma$, being a restriction of the equivalence $L_\ma$, is full and faithful (analogously  for $\mb$).

As is well known every variety equivalent to a variety $_R\Mod$  is necessarily of the form $_S\Mod$ since module categories can be characterized as those varieties  which are Abelian categories and this property is preserved by equivalences. 
(Algebraically one may argue as follows: since every (varietal) generator $G$ in $_R\Mod$ is a  $R$-$S$-bimodule for $S=\mathit{End}(G)$, the endomorphism ring of $G$ (see e.g. \cite[Theorem 17.8]{AF}), every variety $\BV$, which is equivalent to $_R\Mod$ is concretely equivalent to $_S\Mod$.) 
Hence,  one can restrict the above to the bicategory $\mathit{RING}$ of rings with bimodules as 1-cells  and so obtains the classical Morita theory as an immediate consequence.    The varietal generators in  $_R\Mod$  are the so-called \em progenerators\em.

\item[\sf General and special linear groups.] 
Assigning to a commutative unital ring $R$ the monoid $M(n,R)$ of $n\times n$-matrices over $R$, the general linear group $GL(n,R)$ or the special linear group $SL(n,R)$ defines  functors $M_n\colon \Ring\xra{ }\Mon$,  $GL_n\colon \Ring\xra{ }\Grp$, and   $SL_n\colon \Ring\lra\Grp$, respectively. These functors are $\Mon$- and   $\Grp$-representable, respectively and, hence, all have a left adjoint. They do not arise from the canonical constructions.

Since $\N$ is the free monoid over a singleton $\{\star\}$ it suffices  to find a $M_n$-universal monoid homomorphism 
$\N\xra{ u} M(n, A_n)$. Choose  $A_n$ to be $\Z[X_{i,j}; 1\leq i,j\leq n]$, the free commutative ring $Fn^2$ over $n^2$ and  $u$ the monoid morphism mapping $1$ to  the $n\times n$ identity matrix $E_n=(\delta_{ij})$ over $A_n$. For every  $M\in M(n,R)$ we have the unique ring homomorphism $\psi\colon \Z[X_{i,j}]\xra{ }R$  with $\delta_{ij}\mapsto m_{ij}$. Then $M_n(\psi) = M$ and $\psi$ is unique with this property. Thus, $M_n$ is $\Mon$-representable with representing object $Fn^2$. 

Since $GL_n = \Ring\xra{ M_n}\Mon\xra{ (-)^\times}\Grp$ and $(-)^\times$ is right adjoint to the forgetful functor $\Grp\xra{ }\Mon$ the functor $GL_n$ is $\Grp$-representable with representing object $Fn^2$.

Concerning $SL_n$ one proceeds analogously.  One needs to find a $SL_n$-universal group homomorphism 
$\Z\xra{ u} SL(n, A_n)$.  Consider the coequalizer $F(n^2)\xra{ q}A_n$ of the  ring homomorphisms $\det^\sharp, 1^\sharp\colon F1\xra{ }Fn^2$   given by $\det$ and the constant $1$,  considered as  $n^2$-ary derived operations (hence elements of $Fn^2$) in the theory of commutative rings.  It then is easy to see that the group homomorphism $u$ mapping $1\in\Z$ to the $n\times n$ identity matrix $(\delta_{ij})$ over $A_n$ does the job.

In fact Theorem \ref{theo:pres}  more generally implies that every affine group scheme, that is, every  $\Grp$-representable functor on $_c\Alg_R$, for a commutative ring R, has a left adjoint.

\item[\sf  The primitive Hopf structure on polynomial algebras.] 
It is well known  that, for every commutative ring $R$, the polynomial algebras in $_c\mathit{Alg}_R$, {that is,  the free algebras   in $ _c\mathit{Alg}_R$,} carry the structure of an $R$-Hopf algebra in which the variables, that is, the free generators, are primitive elements (see e.g. \cite[p. 92]{Abe}). Recalling that the category  ${_{bi}{\Hopf_{\!R}}}$   of  bicommutative Hopf algebras over $R$ is essentially nothing but $\Coalg(\Ab,{_c\mathit{Alg}_R})$  (see Section \ref{sec:ex}), one gets the  functor $\Set\xra{V_\Phi }{_{bi}{\Hopf_{\!R}}}$  determined by the algebraic functor ${_c\mathit{Alg}_R}\xra{ \Phi^\ast}\Ab$, assigning to a set the respective polynomial algebra.

Thus, item 2 of  Remarks \ref{rem:GS} explains the result mentioned at the beginning of this example, since  $G_\Phi H$ is the set of primitive elements of a bicommutative Hopf algebra $H$  (see Example~\ref{grp-l}).

\item[\sf Adjoint monads induced by bialgebras.]
Any $\T$-bialgebra $\ma$, that is, any 1-cell in $\mathit{LAWV}$ may, as  in any bicategory, allow for the structure of a monad $\Bbb{A}=(\ma, Y^\T\xra{ \eta}\ma, \ma\odot\ma\xra{\mu }\ma)$. Such monad is the same as a monoid in the monoidal category
 $([\T,\T],\odot,Y^\T)$  obtaining its monoidal structure from the monoidal structure of $(\Ladj(\Alg\T,\Alg\T),\circ,id)$ by the equivalence of Theorem \ref{theo:pres}. 
 The monad $\Bbb{A}$, considered as a monad in $(\Ladj(\Alg\T,\Alg\T),\circ,id)$, is a usual monad $(L_\ma, L_\ma\circ L_\ma \xra{\mu}L_\ma, id_{\Alg\T}\xra{\eta}L_\ma)$ on $\Alg\T$ (not to be confused with the monad $\Bbb{A}$ on $\Set$ of Section \ref{sec:monad})  
 and so determines its Eilenberg-Moore category $(\Alg\T)^\Bbb{A}$. 
 It is not difficult to see that the composition of forgetful functors $(\Alg\T)^\Bbb{A}\xra{ }\Alg\T\xra{ }\Set$ is a finitary monadic functor such that one can identify the category  $(\Alg\T)^\Bbb{A}$ with $\Alg\T_\Bbb{A}$ for some Lawvere theory $\T_\Bbb{A}$ (see e.g. \cite[A.21]{ARV}). The forgetful functor of $\Alg\T_\Bbb{A}\xra{ }\Set$ factors as $\Alg\T_\Bbb{A}\simeq (\Alg\T)^\Bbb{A}\xra{ }\Alg\T\xra{|-| }\Set$,  
where the functor $\Alg\T_\Bbb{A}\xra{ U_\mathbb{A}}\Alg\T$ is  finitary monadic; since  $U_\mathbb{A}$ commutes with the forgetful functors it is an algebraic functor and so determines a theory morphism $\T\xra{\Phi }\T_\Bbb{A}$.  Categories of the form $\Alg\T_\Bbb{A}$   first appeared in \cite{TW} with $\T$  the category of commutative rings; here also applications of this construction are given. 

Recall from \cite{EM}:  If $F\dashv U$ be an adjunction from $\BC$ to $\BD$, whose  induced monad on $\BC$ is $\Bbb{A}$,  then the following conditions are equivalent,
\begin{enumerate}
\item There exists an adjunction $U\dashv G$, inducing a comonad  $\Bbb{C}$ on $\BC$.
\item There exists an adjunction $T:= UF\dashv C:=UG$.  
\end{enumerate}
and imply that the (co)Eilenberg-Moore categories $\BC^\Bbb{A}$ and $\BC_\Bbb{C}$ coincide up to a concrete isomorphism over $\BC$ and, hence, $U$ is monadic and comonadic.

Applying this to the above, with $\BC=\Alg\T$, $U_\mathbb{A}=\Phi^\ast$, and $F=\Phi_\ast$ we see, since $L_\ma = U_\mathbb{A}F = \Phi^\ast\circ\Phi_\ast$ has the right adjoint $R_\ma$,  that the functor $U_\mathbb{A}=\Phi^\ast$ is monadic and comonadic.
We so obtain the following characterization of theory morphisms, whose induced algebraic functor has a right adjoint where
the equivalence of 1. and 2. is well known. 

Equivalent are, for any theory morphisms $\T\xra{\Phi }\BS$,
\begin{enumerate}
\item $\Phi^\ast$ has a right adjoint $R$.
\item $\Phi^\ast$ is monadic and comonadic.
\item  The $\T$-bialgebra $\ma$ corresponding to the left adjoint functor $L:= \Phi^\ast\circ \Phi_\ast$ (equivalently, representing the right adjoint $R\circ\Phi^\ast$) carries a monad $\mathbb{A}$ such that  $\BS \simeq \T_\mathbb{A}$ and $\Phi^\ast = U_\mathbb{A}$.
\end{enumerate}
The following equivalence, for  any $\T$-bialgebra $\ma$,  
characterizes the monad structures on left adjoint endofunctors of  varieties in terms of their right adjoints.
$$ L_\ma \text{\ carries the structure of a monad \ }\mathbb{A} \iff 
\text{\ The functor } R_\ma \text{\ has a right adjoint}$$
We illustrate the above by the following simple example: Consider the functors $\Mon\xra{ (-)^\times}\Grp\xra{\Lambda^\ast} \Mon$ of Diagram \eqref{diag:12}. Here $\Lambda^\ast$ has $(-)^\times$ as its right adjoint. $(-)^\times$ is the $\Grp$-lift of the functor $\Mon((\Z,+,0),-)$ and $\Lambda^\ast\circ (-)^\times$ is the $\Mon$-lift of  $\Mon((\Z,+,0),-)$. 

$\Grp$ not only is concretely equivalent to the category $\Mon^\Bbb{T}$ for the monad $\Bbb{T}$ given by the adjunction $\Lambda_\ast\dashv\Lambda^\ast$, but also to $\Mon_\Bbb{A}$ for the comonad $\Bbb{A}$ given by the adjunction $(-)^\times \vdash \Lambda^\ast$: a monoid $M$ carries a (unique) group structure iff  $M\simeq\Lambda^\ast (M^\times)$ 
that is, in other words, groups are the coalgebras of the comonad $\Bbb{A}$, since the counit of  the adjunction $(-)^\times \vdash \Lambda^\ast$ is an isomorphism. 
\end{description}

\subsection{Coalgebras in commutative varieties}\label{sec:com}
A Lawvere theory $\T$ is called \em commutative\em,  if for all $\sigma\in\T(n,1), \tau\in \T(m,1)$ the following diagram commutes.  
\begin{equation*}
\begin{aligned}
\xymatrix@=1.5em{
 {n\times m}    \ar[r]^{\ \ \  \tau^n  }\ar[d]_{ \sigma^m }&{n} \ar[d]^{\sigma }  \\
 m   \ar[r]_{ \tau}          &    1 
}
\end{aligned}
\end{equation*}
Varieties $\BV$ whose Lawvere theory $\TV$ is commutative are known under various names as, e.g., \em entropic varieties\em, \em distribute varieties \em or \em commutative varieties \em (see \cite{Davey} for the respective references).
{Though the first option seems to be the more popular choice amongst universal algebraists, we will use the second one because of its match with the notion of a commutative theory.}
We notice that for a commutative theory $\T$ the categories $\Alg(\T,\Alg\T)$ and $\Alg\T$ are equivalent.

Since limits in functor categories are computed point wise one immediately gets from this definition, that for  any commutative Lawvere theory $\T$ and  for every internal $\T$-algebra $\ma$ in  a  category $\BC$ with finite products the $\ma$ interpretation $\ma \tau$ of any $m$-ary operation $\tau\in\T(m,1)$ is a morphism $\ma^m\xra{ }\ma$ in $\Alg(\T,\BC)$, that is, for every $\sigma\in\T(n,1)$ and every $n\in\N$ the following diagram commutes in $\BV$, where $A:=\ma 1$.   \vspace{-1em}
\begin{equation*}
\begin{aligned}
\xymatrix@=2em{
 A^{n\times m}    \ar[r]^{\ \ (\tau^\ma)^n  }\ar[d]_{\sigma^{(\ma^m)} =(\sigma^\ma)^m }&{A^n } \ar[d]^{\sigma^\ma  }  \\
 A^m   \ar[r]_{ \tau^\ma}          &    A 
}
\end{aligned}
\end{equation*} 

By dualization we obtain
\begin{Proposition}\label{prop:com}
 {Let $\T$ be a commutative  Lawvere theory. Then for every variety $\BV$} the category $\Coalg(\T,\BV)$ is \em cocommutative\em, that is, for every  $\T$-algebra $\T\xra{ \ma}\BV^\op$ and every $\tau\in\T(m,1)$ the $\BV$-homomorphism $\ma\tau = A\xra{ \tau_\ma}m\cdot A$ is a morphism in $\Coalg(\T,\BV)$.
\end{Proposition}

\subsubsection{Canonical coalgebras}\label{ssec:com}

Since for a commutative variety $\BV$ we have the equivalence $\Alg(\T_\BV,\BV)\simeq \BV$, the functor $T$ of Proposition \ref{prop:can} is an equivalence and, hence the functor $T$ is a coreflective embedding. Again, in this situation a simpler description of this fact can be given as follows.

In a commutative  variety $\BV$ the hom-sets $\BV(A,B)$ form subalgebras $[A,B]$ of the products $B^A$. 
 $\BV$ then is a monoidal closed category with  internal hom-functor $[-,-]$ so defined  (see \cite{Davey}, \cite{Bor}). 
 In particular, for each $\BV$-algebra $A$ one has an adjunction $A\ot -\dashv [A,-]$ on $\BV $ and, thus, a \em canonical \em $\T_\BV$-coalgebra $NA$ with underlying $\BV$-algebra $A$.  Equivalently, $NA$ is the $\TV$-coalgebra in $\BV$ 
with corresponding $\T_\BV$-representable functor $R_{NA} = \BV\xra{ [A,-]}\BV$.
The co-operation $A\xra{ \sigma_{NA}}n\cdot A$ of $NA$   corresponding to  $\sigma\in\TV(n,1)$ is (use  Equation \eqref{eqn:coopgen} and the fact that operations in $\BV$ are homomorphisms and Theorem  \ref{theo:pres}, respectively) 
  \vspace{-0.5em}
\begin{equation*}
A\xra{ \sigma_{NA}} n\cdot A =(A\simeq A\otimes F_\BV 1\xra{A\otimes\sigma } A\otimes F_\BV n \simeq n\cdot A)
= A \xra{ \langle\nu_1,\ldots,\nu_n\rangle} (n\cdot A)^n\xra{  \sigma^{n\cdot A}} n\cdot A
  \vspace{-0.5em}
\end{equation*}
This defines a full embedding $N\colon \BV\xra{ }\Coalg(\T_\BV,\BV) \text{\ with \ } |\!|-|\!|\circ N = id_\BV$. 
Though $N$ occasionally is an equivalence, that is, up to isomorphism $\BV$ only admits canonical $\TV$-coalgebras, this is not the case general (see Examples \ref{ex:example1} below). 
Slightly more general we have the following results.
\begin{Theorem}\label{thm:corefl1}
Let $\Phi\colon \T\xra{ }\TV$ be a theory morphism into a commutative theory $\TV$. Then the following hold, with notation as above:
\begin{enumerate}
\item The functor $N_\Phi = \BV\xra{N}\Coalg(\TV,\BV)\xra{{^{\!\BV}\Phi}}\Coalg(\T,\BV)$  is  full and faithful.
\item $N_\Phi( F_\BV X) = V_\Phi(X)$ and  $|\!|-|\!|\circ N_\Phi = \id_\BV$.
\item The assignment $\Phi\mapsto N_\Phi$ defines  an essentially  bijective correspondence between morphisms of Lawvere theories $\Phi\colon\T\xra{ }\TV$ and  functors   $S\colon\BV\xra{ }\Coalg(\T,\BV)$ with $|\!|-|\!|\circ S = \id_\BV$.
\item The $\T$-representable functor corresponding to $N_\Phi A$  is $R_{N_\Phi A} =\BV\xra{[A,-]}\BV\xra{\Phi^\ast}\Alg\T$, for any $\BV$-algebra $A$. The co-operation $\sigma_{N_{\Phi}A}$ corresponding to $\sigma\in \mathcal{T}(n,1)$ is explicitly given by $\sigma_{N_{\Phi}A}=\Phi(\sigma)^{n\cdot A}\circ \langle\nu_1,\cdots,\nu_n\rangle.$
\end{enumerate}
 \end{Theorem}
 
\begin{Proof}
Considering coalgebras as representable functors $N_\Phi$ acts as $A\mapsto \Phi^\ast\circ [A,-]$ and
 $N_\Phi f\colon  \Phi^\ast [B,-]\Rightarrow\Phi^\ast [A,-]$ is the natural transformation with components $\Phi^\ast [B,C]\xrightarrow{\Phi^\ast [f,C]}\Phi^\ast [A,C]$, for any $\BV$-morphism $A\xra{ f}B$. Now the natural transformations $|N_\Phi f|\colon |\Phi^\ast [B,-]|=\BV(B,-)\Rightarrow \BV(A,-)=|\Phi^\ast [A,-]|$ correspond one-to-one to $\BV$-morphisms $A\xra{ f}B$, which proves~1.
The first identity of item 2 follows from Equation \eqref{eqn:coopgen} since $\BV$ is commutative, hence all operations are homomorphisms; the second one is obvious, as is item 4. 

For every functor  $S\colon\BV\xra{ }\Coalg(\T,\BV)$ with $|\!|-|\!|\circ S = id_\BV$ there exists by Proposition \ref{corr:algfu} the algebraic functor  $R_{SF_\BV 1}\colon \BV\xra{ }\Alg\T$. 
Let ${\Phi}_S\colon\T\xra{ }\TV$ be the  morphism of Lawvere theories with  ${\Phi_S}^\ast = {R_{S F_\BV 1}}$. The assignment $S\mapsto \Phi_S$ is  essentially injective: indeed, for each $\BV$-homomorphism $F_\BV 1\xra{ f}A$   the $\T$-coalgebra morphism $Sf$ has $f$ as its underlying morphism in $\BV$ and this satisfies,  for each each $\sigma \in \T(n,1)$, the equation $F_\BV1 \xra{ f} A \xra{ \sigma_{SA}}n\cdot A = F_\BV1 \xra{ \sigma_{SF_\BV 1}} A \xra{n\cdot f }n\cdot A$.
This proves that $S$ is determined by  $S F_\BV 1$ since the family $(F_\BV 1\xra{ f}A)_{f\in\BV(F_\BV 1, A)}$ is jointly epimorphic.  
Items 2. and 3. of  Theorem \ref{thm:corefl1} imply $ R_{N_\Phi F_\BV 1} \simeq \Phi^\ast$ for each $\Phi$. In other words,  the assignment $S\mapsto \Phi_S$ is essentially bijective.
\end{Proof}
\begin{Theorem}\label{thm:corefl2}
Let $\Phi\colon \T\xra{ }\TV$ be a theory morphism into a commutative theory $\TV$. Then the following hold. 
\begin{enumerate}
\item For every $\T$-coalgebra $\ma$ in $\BV$ the set $G_\Phi(\ma)$ is the underlying set of a  $\BV$-subalgebra of $|\!|\ma|\!|$, such that the functor $G_\Phi$ factors as $\Coalg(\T,\BV)\xra{ \bar{G}_\Phi}\BV\xra{ |-|}\Set$ and the natural transformation $e$ lifts to a natural transformation $\bar{e}\colon \bar{G}_\Phi\Rightarrow |\!|-|\!|$. 
The  lifted functor $\bar{G}_\Phi$  satisfies  the equations $\bar{G}_\Phi N_\Phi = \id_{\BV}$ and  $\bar{G}_\Phi V_\Phi(X) = F_\BV X$ for each set $X$.
\item $N_\Phi$ {is left} adjoint {to} $\bar{G}_\Phi$; hence $\BV$ is (equivalent to) a full coreflective subcategory of $\Coalg(\T,\BV)$.
\end{enumerate}
 \end{Theorem}
 \begin{Proof}
Concerning item 1 observe first that $G_\Phi(\ma)$ is the underlying set of a $\BV$-subalgebra of $A$   by Equation \eqref{eq:G}, since  the maps $\Phi(\tau)^{n\cdot A}$ are homomorphisms in  $\BV$.  In particular, $G_\Phi$ factors over $\BV$. The lifted functor $\Coalg(\T,\BV)\xra{\bar{G}_\Phi}\BV$  obviously  satisfies the equation  $\bar{G}_\Phi N_\Phi = id_{\BV}$. 
Finally, item 1 implies, again by  the lifting theorem of adjunctions,  that the functor $\bar{G}_\Phi$ has a left adjoint $L$, which  then satisfies  the equation $L( F_\BV X) = V_\Phi(X)$. Now $N_\Phi(F_\BV 1) = L(F_\BV 1)$ follows by item 2 of Theorem \ref{thm:corefl1}, and this implies $L=N_\Phi$ by item 4 of Theorem \ref{thm:corefl1}.
\end{Proof}

\begin{Remark}\label{rem:isos}\rm 
For any theory morphism $\Phi\colon \T\xra{ }\TV$   into a commutative theory $\TV$ the following are equivalent:
\begin{enumerate}
\item $N_\Phi$ is an equivalence with equivalence inverse $\bar{G}_\Phi$. 
\item {${G}_\Phi\simeq |-|\circ |\!|-|\!|$.}
\item $N_\Phi \simeq C_{\T,\BV}$. 
\end{enumerate}
\end{Remark}

The following examples show that the functors $N_\Phi$ may or may  not be equivalences.
\begin{Examples}\label{ex:example1}
 \begin{enumerate}
 \item\label{free} 
 If $\BV$ is a commutative  variety all of whose algebras are free, then $N$ is an equivalence as is easily seen. 
  By \cite{Givant} there are essentially four such varieties: $\Set$, the variety of pointed sets $\Set_\ast$, $\Vect$ and ${\mathit{Aff}_k}$, the varieties of vector spaces  and affine spaces, respectively, over a field $k$. 
Thus, if $\BV$ is any of these varieties, then  every $\T_\BV$-coalgebra in $\BV$ is of the form $\BV(X)$ as it is the case for the variety of groups (see Remark \ref{rems:Kan}). 

 \item\label{ex:examples} 
 By the last example in Section \ref{sec:ex} there are the   equivalences 
$\Coalg(\Ab,\Mod_R)\simeq \Mod_R\simeq\Coalg(\Grp,\Mod_R)$ and  $\Coalg({_c\Mon},\Mod_R)\simeq \Mod_R\simeq\Coalg(\Mon,\Mod_R)$,  which we will use below. 
The  forgetful functors coincide with $\bar{G}_\Phi$ and are equivalences, concrete over  $_c\Mon$ and $\Ab$, respectively, having the  functors  $N_\Phi$ as their inverses.
   
 \item Let $\BV={_c\mathit{Sem}}$ be the variety of commutative semigroups, and  $\Phi^\ast\colon \Ab\to{_c\mathit{Sem}}$ the forgetful functor. Note that $_c\mathit{Sem}$ is commutative but fails to be semi-additive. 
 
 The coproduct ${S}+{T}$ in $\BV$ 
  is  {given} as $({S}^0\times {T}^0)\setminus\{\, (0,0)\,\}$ with the component-wise addition, where ${S}^0$ (respectively  ${T}^0$) denotes the unitarization of ${S}$ (resp. ${T}$), and with coproduct injections $\nu_1(a):=(a,0)$, $a\in S$, and $\nu_2(b):=(0,b)$, $b\in T$. 
  Then ${N}({S})$ is an internal $\mathcal{T}_{_c\mathit{Sem}}$-coalgebra in $_c\mathit{Sem}$ with comultiplication $\mu_{N\mathsf{S}}=\mu^{2\cdot {S}}\circ \langle\nu_1,\nu_2\rangle$, i.e., $\mu_{N{S}}(a)=(a,0)+(0,a)=(a,a)$, $a\in S$. 
 \begin{enumerate}
 \item\label{ex:pt4}  
Let ${S}$ be a commutative semigroup that contains an idempotent element $e$ (e.g.  the underlying semigroup of a commutative monoid).
Let $E_e\colon {S}\to {S}+{S}$ be the semigroup homomorphism given by $E_e(a):=(e,e)$, $a\in S$. This endows $\mathsf{S}$ with a structure of an internal $\mathcal{T}_{_c\mathit{Sem}}$-coalgebra in $_c\mathit{Sem}$ not isomorphic to that  obtained by applying $N$ when $S$ is not reduced to $e$. This shows that $N$ cannot be an equivalence between $_c\mathit{Sem}$ and $\Coalg(\mathcal{T}_{_c\mathit{Sem}},{_c\mathit{Sem}})$.
 \item \label{ex:pt3}
  The comultiplication in ${N}_{\Phi}({A})$, for an abelian group ${A}$, thus, is  given by $a\mapsto (a,a)$. But ${A}$ carries also a trivial internal $\mathcal{T}_{\mathit{cSem}}$-coalgebra structure, namely, $Z_{{A}}(a):=(0,0)$. This clearly implies that ${N}_{\Phi}$ is not an equivalence of the categories  $\Ab$ and $\Coalg(\mathcal{T}_{_c\mathit{Sem}},\Ab)$.
\end{enumerate}
 \end{enumerate}
 \end{Examples}

 \subsubsection{Further applications}\label{sec:more}
 
 \begin{Lemma}\label{lem:factorization}
Let $\Phi= \T_\BW\xra{\Xi }\T_\BU \xra{ \Sigma} \TV$ be a factorization of a morphism of Lawvere theories over a commutative theory $\T_\BU$.
Then the following diagram commutes   \vspace{-0.5em}
\begin{equation*}\label{diag:factorizations}
\xymatrix@=2em{
\Coalg(\mathcal{T}_{\mathcal{W}},\mathcal{V})\ar[dr]_{G_{\Phi}} \ar[r]^{{_{\BW}\Sigma_\ast}} & \Coalg(\mathcal{T}_{\mathcal{W}},\mathcal{U})\ar[r]^{\ \ \ \ \ \bar{G}_{\Xi}}\ar[d]_{G_{\Xi}} & \mathcal{U}\ar[dl]^{|-|}\ar[d]^{\Xi^{\ast}} \\
 & \Set & \mathcal{W} \ar[l]^{|-|} & 
}\vspace{-1em}
\end{equation*}
and, moreover, 
\begin{enumerate}
\item $\bar{G}_\Xi ({_{\BW}\!\Sigma_\ast} \ma)$ is a $\mathcal{U}$-subalgebra of $\bar{G}_{\Xi}({_{\BW}\!\Sigma_\ast}(C_{\mathcal{T}_{\mathcal{W}},\mathcal{V}}(|\!|\ma|\!|))) =  \bar{G}_{\Xi}(C_{\mathcal{T}_{\mathcal{W}},\mathcal{U}}(\Sigma^{\ast}(|\!|\ma|\!|)))$
 for every $\T_\BW$-algebra $\ma$ in $\BV$, 
\item $\bar{G}_\Xi ({_{\BW}\!\Sigma_\ast} \ma)$ is a $\mathcal{U}$-subalgebra of $ \bar{G}_{\Xi}((\Sigma^{\ast}(|\!|\ma|\!|)))$, if 
$\bar{G}_\Xi$ is an isomorphism.
\end{enumerate}
\end{Lemma}

\begin{Proof}
The first statement is obvious (see Lemma   \ref{lem:comp}). Concerning the second one note that, with $\ma\xra{ \eta_\ma} |\!|C_{\mathcal{T}_{\mathcal{W},\BV}}|\!|$  the adjunction unit, the injective map  $||\!|\eta_\ma|\!| |\circ e_\ma$ factors as  $e_{C_{\mathcal{T}_{\mathcal{W}},\mathcal{V}}(|\!|\ma|\!|)}\circ G_\Phi(\eta_\ma)$ by naturality (see Equation \eqref{eq:G}).
This implies that $\bar{G}_\Xi ({_{\BW}\!\Sigma_\ast} (\eta_\ma)$ is a $\BU$-monomorphism. Item~3. follows trivially (see also Remark \ref{rem:isos}).
\end{Proof}

In the discussion of the following applications we use  notations as in Diagram \eqref{diag:12}. 
\begin{description}
\item[{\sf More on primitive elements.}] 
\begin{enumerate}
\item The primitive element functor ${_{c}\Bialg_R}\xra{G_\Phi}\Set$ has by Lemma \ref{lem:factorization} the following factorization where, by Examples \ref{ex:example1}.\ref{ex:examples},   $\bar{G}_\Xi$ is an isomorphism with inverse $N_\Xi$ and, thus, can be identified with the functor $|\!|-|\!|$. 
\vspace{-0.5em}
 $$\Coalg(\Mon,{_c\mathit{Alg}_R})\xrightarrow{{_{\mathcal{M}}\Sigma_\ast}}\Coalg(\Mon,\Mod_R)\xrightarrow{\bar{G}_{\Xi}}\Mod_R\xrightarrow{|-|}\Set.  \vspace{-0.5em}$$
We so obtain the following familiar result.
\em For each commutative Hopf algebra $A$ the set  $G_{\Phi}(A)$ of primitive elements of $A$ is a submodule of the
the carrier algebra $\Sigma^{\ast}(A)$\em. 
\item Recall that each symmetric  algebra, that is, each free commutative $R$-algebra $\Sigma_\ast(M)$ over an $R$-module $M$, carries the structure of a Hopf algebra  $S(M)$, where $\Mod_R\xra{ \Xi^\ast}\Ab$ is the forgetful functor. 
This construction  is given by the functor
\vspace{-0.5em}
$$S=\Mod_R\xra{ N_\Xi}\Coalg(\Ab,\Mod_R)\xra{ ^{\BW\!}\Sigma_\ast}\Coalg(\Ab,{_c\mathit{Alg}_R})\simeq {_{bi}\Hopf_{\!R}}  \vspace{-0.5em}$$ 
 $S$ is left adjoint to $G_\Xi\circ{_{\BW}\Sigma_\ast}$ by  Proposition~\ref{prop:props}; since $N_\Xi$ has $G_\Xi$ as its inverse (see Examples \ref{ex:example1}), this adjunction is essentially the adjunction $^{\BW\!}\Sigma_\ast\dashv {_{\BW}\Sigma_\ast}$. Since the unit of this adjunction is point-wise monomorphic (see the proof of Proposition  \ref{prop:props}), we obtain the familiar result that \em the set of primitive elements of  $S(M)$ contains $M$\em.
\end{enumerate}
\item[\sf More on group-like elements.] Since $\Coalg(\Mon,{_c\Mon})=\Coalg(_c\Mon,{_c\Mon})$ (see  Examples \ref{ex:example1}.\ref{ex:examples}),  the group-like element functor ${_{c}\Bialg_R}\xra{G_\Psi}\Set$ 
  factors by Lemma \ref{lem:factorization}  as  
  \vspace{-0.5em}
  $$\Coalg(\Mon,{_c\mathit{Alg}_R})\xrightarrow{_{\T_{\mathcal{M}}}{_c\Psi_\ast}}\Coalg(\Mon,{_c\Mon}) \xrightarrow{\bar{G}_{_c\mathcal{M}}} {_c\Mon} \xrightarrow{|-|} \Set.   \vspace{-0.5em}$$
 where  the functor $\bar{G}_{_c\mathcal{M}}$ is an isomorphism. We so obtain   the   fact that \em the set of group-like elements of a commutative bialgebra (hence, in particular, of any  commutative Hopf algebra) $A$  is a submonoid of the multiplicative monoid of the underlying algebra of $A$. \em
\end{description}

\end{document}